\documentclass[10pt]{article}
\usepackage[utf8]{inputenc}

\usepackage[top=1in,bottom=1.25in,left=1.25in,right=1.25in]{geometry}

\usepackage[absolute,overlay]{textpos}

\usepackage{hyperref}
\hypersetup{
  colorlinks=true,
  linkcolor=black,
  filecolor=red,      
  urlcolor=blue,
  citecolor=blue
}

\usepackage{graphicx}
\usepackage{ragged2e}
\usepackage{tabto}
\usepackage{float}
\usepackage{multirow}
\usepackage{amsmath}
\usepackage{amssymb}
\usepackage{amsthm}
\usepackage{mathtools}
\usepackage[ruled,vlined]{algorithm2e}
\usepackage{booktabs}

\usepackage{tikz}
\usepackage{pgfplots}
\usetikzlibrary{arrows.meta}
\usepackage{xcolor}
\definecolor{darkred}{RGB}{173,0,43}

\usepackage[nottoc]{tocbibind}

\usepackage[font=footnotesize]{subcaption}
\usepackage[font=footnotesize]{caption}

\usepackage{multirow}
\usepackage{multicol}

\usepackage[capitalise,noabbrev]{cleveref}
\crefname{equation}{}{}
\usepackage{afterpage}
\usepackage{marginnote}

\usepackage{lmodern}

\usepackage[sf,bf,small]{titlesec}

\usepackage{bm}
\usepackage{dsfont}

\newcommand\vct[1]{\bm{\mathsf{#1}}}
\newcommand\mtx[1]{\bm{\mathsf{#1}}}

\newcommand{\lp}{\left(}
\newcommand{\rp}{\right)}

\usepackage[sort,compress]{cite}

\newcommand\bbR{\mathbb R}

\newcommand\bx{\bm{x}}
\newcommand\bs{\bm{s}}
\newcommand\by{\bm{y}}
\newcommand\bt{\bm{t}}

\newcommand\bn{\bm n}

\newcommand\bC{\mathbb C}

\newcommand\bX{\bm X}

\newcommand\tinit{t_q}

\newcommand\sinit{S_q}

\newtheorem{remark}{\sffamily Remark}
\newtheorem{definition}{\sffamily Definition}

\newcommand{\cI}{\mathcal I}
\newcommand{\cD}{\mathcal D}
\newcommand{\cS}{\mathcal S}

\newcommand{\cO}{\mathcal O}

\newcommand{\cF}{\mathcal F}

\newcommand{\cK}{\mathcal K}

\newcommand{\sfB}{\mathsf B}
\newcommand{\sfN}{\mathsf N}
\newcommand{\sfF}{\mathsf F}
\newcommand{\sfQ}{\mathsf Q}
\newcommand{\sfP}{\mathsf P}
\newcommand{\sfI}{\mathsf I}
\newcommand{\sfJ}{\mathsf J}
\newcommand{\sfK}{\mathsf K}
\newcommand{\sfS}{\mathsf S}
\newcommand{\sfR}{\mathsf R}

\newcommand{\ext}{\textrm{ext}}

\newcommand{\note}{\color{black}}

\newcommand{\Npatches}[0]{N_{\textrm{patches}}}
\newcommand{\ucomp}[0]{u_{\textrm{comp}}}
\newcommand{\ppw}{{\rm ppw}}
\newcommand{\ppwmin}{{\rm ppw}_{{\rm min}}}

\renewcommand{\phi}{\varphi}

\numberwithin{equation}{section}

\begin{document}

\begin{titlepage}

  \raggedleft
  {\sffamily \bfseries Status: {\color{red} arXiv pre-print}}
  
  \hrulefill
  
  \raggedright

  \begin{textblock*}{\linewidth}(1.25in,2in) 
    {\LARGE \sffamily \bfseries FMM-LU: A fast direct solver for multiscale boundary  \\
      \vspace{.25\baselineskip}
      integral equations in three
      dimensions}
  \end{textblock*}

  \vspace{1.5in}
  Daria Sushnikova\footnote{Research supported in part by
  the Office of Naval Research under award
  numbers~\#N00014-17-1-2451 and~\#N00014-18-1-2307.}\\
\small \emph{King Abdullah University of Science and Technology\\
Thuwal 23955-6900, Saudi Arabia}\\
  \texttt{daria.sushnikova@kaust.edu.sa}
  \normalsize

  \vspace{\baselineskip}
  Leslie Greengard\footnote{Research supported in part by 
  the Office of Naval Research under award
    number~\#N00014-18-1-2307.}\\
  \small \emph{Courant Institute, NYU\\
    New York, NY, 10012}\\
  \texttt{greengard@cims.nyu.edu}
  \normalsize

   \vspace{\baselineskip}
   Michael O'Neil\footnote{Research supported in part by
    the Office of Naval Research under award
    numbers~\#N00014-17-1-2451 and~\#N00014-18-1-2307.}\\
  \small \emph{Courant Institute, NYU\\
    New York, NY 10012}\\
  \texttt{oneil@cims.nyu.edu}
  \normalsize
  
  \vspace{\baselineskip}
  Manas Rachh\\
  \small \emph{Center for Computational Mathematics, Flatiron Institute\\
    New York, NY 10010}\\
  \texttt{mrachh@flatironinstitute.org}
  \normalsize

  \begin{textblock*}{\linewidth}(1.25in,7in) 
    \today
  \end{textblock*}

  
\end{titlepage}

\begin{abstract}
  
  We present a fast direct solver for boundary integral equations on complex
  surfaces in three dimensions using an extension of the recently introduced
  recursive strong  skeletonization scheme. For problems that are not highly
  oscillatory, our algorithm computes an $\mtx{LU}$-like hierarchical
  factorization of the dense system matrix, permitting application of the
  inverse in $\mathcal O(n)$ time, where $n$ is the number of unknowns on the surface.
  The factorization itself also scales linearly with the system size, albeit
  with a somewhat larger constant. The scheme is built on a level-restricted
  adaptive octree data structure, and therefore it is compatible with highly
  nonuniform discretizations. Furthermore, the scheme is coupled with high-order
  accurate locally-corrected Nystr\"om quadrature methods to integrate the
  singular and weakly-singular Green's functions used in the integral
  representations. Our method has immediate applications to a variety of
  problems in computational physics. We concentrate here on studying its
  performance in acoustic scattering (governed by the Helmholtz equation) at low
  to moderate frequencies, and provide rigorous justification for compression
  of submatrices via proxy surfaces.\\


\end{abstract}

\tableofcontents

\newpage

\section{Introduction}
\label{sec:intro}

Integral equation formulations lead to powerful methods for the
solution of boundary value problems governed by the partial
differential equations (PDEs) of classical mathematical physics.  The
corresponding free-space Green's functions are well-known, and in the
absence of source terms, boundary integral formulations are restricted
to the surface of the domain, thereby reducing the dimensionality of
the problem.  There has been a resurgence of interest in these methods
over the past few decades because of the availability of fast
algorithms for applying the dense matrices which arise after
discretization.  These include hierarchical schemes such as fast
multipole methods (FMMs), panel clustering methods, $\mathcal H$-matrix methods
and multigrid variants, as well as FFT-based schemes such as the
method of local corrections, pre-corrected FFT methods, etc.  The
literature on such methods is vast and we point only to a few review
articles and monographs
\cite{Borm-h-2003,greengard-1997,shortcourse,martinsson-book,liubook}.
Assuming $n$ is the number of degrees of freedom used in sampling the
surface and the corresponding charge and/or dipole densities, and
$\mtx{A}$ is the $n \times n$ system matrix obtained after the
application of a suitable quadrature rule to the chosen integral
representation, these algorithms permit $\mtx{A}$ to be applied to a vector
in $\cO(n)$ or $\cO(n \log n)$ time. When the linear system is
well-conditioned, this generally allows for the rapid iterative solution of
large-scale problems in essentially optimal time.

Despite these advances, there are several tasks where iterative solvers are not
satisfactory. The obvious case is in problems which they themselves are
ill-conditioned (such as scattering \emph{near resonance}, i.e. from a
cavity-like object with nearly-resonant modes). An equally important area where
direct methods are preferred is in inverse problems, or in any other setting
where one needs to solve the same system matrix with multiple right-hand sides.
Finally, fast direct solvers are very useful when exploring low-rank
perturbations of the geometry and, hence, the system matrix. Updating the
solution in such cases requires only a few applications of $\mtx{A}^{-1}$ or a
fast update of the inverse itself \cite{greengard-2009,minden-2016}.

In the last few years, several algorithmic ideas have emerged which permit the
construction of a compressed approximation of $\mtx{A}^{-1}$ at a cost of the
order $\cO(n)$ or $\cO(n \log^p n)$, for modest $p$. In this paper, we describe
such a scheme, which we refer to as the {\em FMM-LU method}. It uses FMM-type
hierarchical compression strategies to rapidly compute an
$\mtx{LU}$-factorization of the system matrix, building on the algorithmic
framework of the recently introduced \emph{recursive strong skeletonization}
procedure~\cite{minden-str_scel-2017}. In this manuscript, we apply the method
to adaptive, multiscale surface discretizations of boundary integral equations
coupled to high-order accurate quadratures; this leads to efficient,
high-fidelity solvers for geometrically intricate models. In this work, we will
concentrate on boundary value problems for the Helmholtz equation at low to
moderate frequencies, but the scheme is equally applicable to many other
families of boundary value problems, e.g. Laplace, Stokes, Maxwell, etc., In
each case, however. there are many application-specific details that need to be
addressed.

\subsection*{Related work}

Fast direct solvers for integral equations have their origin in the solution of
two-point boundary value problems in one dimension
\cite{gr-bvp,sr-bvp,ambikasaran2013n,ChMing-dir_hss-2006,ChLi-hss-2007}. There
is an extensive literature in this area and the methods are, by now, quite
mature. In higher dimensions, a major distinction in fast algorithms for both
hierarchically compressible matrices and their inverses concerns which blocks
are left intact and which are subject to compression. For the sake of
simplicity, let us assume an octree data structure has been imposed on a surface
discretization and refined uniformly. Let us also assume that the unknowns are
ordered so that the points in each (non-empty) leaf node, denoted by $B_i$, are
contiguous. If all block submatrices $\mtx{A}_{ij}$ corresponding to interactions
between leaf nodes~$B_{i}$ and~$B_{j}$, with $i \neq j$, are (hierarchically)
compressed, using the method of~\cite{cheng_2005,liberty2007randomized} leads to
the recursive skeletonization methods
of~\cite{gillman-dir_hss-2012,ho-2012,martinsson-2005,kong-hodlrdir-2011,Mar-hss-2011}.

These methods are, more or less, optimal for boundary integral equations in two
dimensions but \emph{not} in three dimensions. In three dimensions, the
interaction between neighboring leaf nodes leads to relatively high-rank block
matrices. In the literature, this is referred to as \emph{weak admissibility} or
\emph{weak skeletonization}, and the associated matrix factorizations take the
form of HSS or HODLR formats~\cite{martinsson-hss_integr-2005,martinsson-book,
  hackbusch-h-1999,hackbusch-h-2000,Borm-h-2003}. To overcome this obstacle,
borrowing from the language of fast multipole methods, one can instead choose to
compress only those block matrices~$\mtx{A}_{ij}$ corresponding to {\em
  well-separated} interactions, which are known {\em a priori} to be low-rank to
any fixed precision (based on an analysis of the underlying PDE or integral
equation). Well-separated here means that leaf nodes~$B_i$ and~$B_j$ are
separated by a box of the same size. In the literature, this is referred to as
\emph{strong admissibility}, \emph{strong skeletonization}, or
\emph{mosaic-skeletonization}, and the resulting matrix factorizations take the
form of~\emph{${\cal H}$-matrix} or~\emph{${\cal H}^2$-matrix}
compression~\cite{ambikasaran-2014,Borm-h2-2010,hackbusch-h2-2000,
  Jiao-h2-ce-2017,minden-str_scel-2017,coulier-ifmm_prec-2017}.

Fast solvers for boundary integral equations in three dimensions using weak
skeletonization without additional improvements have led to solvers whose
compression/factorization costs are of the order $\cO(n^{3/2})$, with subsequent
applications of the inverse requiring only~$\cO(n \log n)$ work. Fortunately,
the implicit constants associated with these asymptotics are relatively
small~\cite{ho-2012}. A variety of improvements have been developed to reduce
the factorization costs by introducing auxiliary variables in a
geometrically-aware fashion~\cite{ho-2014,CorMar-dir_hss-2014}.

Strong skeletonization-based schemes are more complicated to handle in terms of
data structures and linear algebraic manipulation, but have significant
advantages. First, a large body of work by Hackbusch and collaborators has led
to a rigorous algebraic theory for hierarchically compressible
${\cal H}^2$-matrices and fast solvers with linear or quasilinear complexity
\cite{Beb-hlu-2005,borm-h2lu-2013,Hakb-h2-lib}. The constants implicit in the
asymptotic scaling with this approach, however, are quite large. To improve
performance, the \emph{inverse fast multipole method} (IFMM) was introduced by
Ambikasaran, Coulier, Darve, and Pouransari
\cite{ambikasaran-2014,coulier-ifmm_prec-2017}, and the strong skeletonization
method was introduced by Minden, Ho, Damle and Ying~\cite{minden-str_scel-2017},
the latter of which we will largely follow here. Recently, a strong
skeletonization-based fast direct solver was shown to be useful in computing
electromagnetic scattering solutions from dielectic
objects~\cite{jiang2022skel}, albeit using low-order discretizations.

Finally, we should note that an early fast solver for {\em volume} integral
equations in two dimensions was presented in~\cite{chen-direct}, and that related
fast solvers have been developed using direct discretization of the PDE rather
than integral equations \cite{XiaSh-hss-2009,hps-direct,solovyev-hss-2014}.
Closely related to our work here is that of \cite{sushnikova-ce-2018} which
makes use of strong admissibility in the sparse matrix setting.

{\note
\subsection*{Contribution}

The algorithm described in the present paper has novel contributions building on
top of the foundational work in~\cite{minden-str_scel-2017}. Each of these
contributions was specifically developed in order to accurately and robustly
solve a particular boundary value problem; in this manuscript we address the
exterior Helmholtz Dirichlet scattering problem. Other problems, such as the
exterior Neumann problem, would require additional modifications, namely to the
method of compression via proxy surfaces. The concepts introduced
in~\cite{minden-str_scel-2017} laid the algorthmic framework for the fast direct
integral equation solver of this work, similar to how standard $n$-body FMM
codes~\cite{wideband3d} form the algorithmic framework for FMM-accelerated
integral equation solvers~\cite{greengard2021fast}. Coupling such algorithmic
frameworks with high-order quadratures and multiscale geometries is non-trivial,
and there are various considerations that must be taken into account.

Here we highlight the novel contributions of this paper which allow for solving
the boundary value problem via an integral equation formulation:
\begin{itemize}
  \item \emph{Multiscale geometries with adaptive discretizations}: The FMM-LU
        algorithm of this work is designed specifically to compress and invert
        matrices arising from adaptive discretizations of boundary integral
        equations along surfaces with multiscale features. We make use of an
        adaptive octree data structure to keep track of analytic estimates
        concerning the maximal rank of well-separated blocks which could be at
        different levels in the tree hierarchy. The resulting hierarchical
        elimination scheme is affected by the adaptive discretization along the
        surface: the adaptive data structure affects the manipulation of Schur
        complements in the ~$\mtx{LU}$-factorization.

  \item \emph{High-order local quadrature corrections}: In order to accurately
        approximate the weakly-singular integral operators for various boundary
        integral equations, it is necessary to use high-order quadrature
        corrections. The resulting linear system is therefore not merely a
        collection of kernel evaluations: entries corresponding to near field
        interactions are modified based on these quadrature corrections. In
        order to couple such a linear system with a fast direct solver, it is
        required to develop quadrature machinery which permits {\em on-the-fly}
        extraction of near field matrix elements when targets are on-or-close-to
        a surface triangle where accurate approximation of a layer potential
        requires some care. The key ingredient here is the use of generalized
        Gaussian quadrature~\cite{bremer,bremer_2013} for self-interactions on
        surface panels and precomputed hierarchical interpolation matrices to
        accelerate adaptive integration in the near field.

  \item \emph{Proxy surface compression theory}: In order to accelerate the
        compression (or factorization) step of a fast direct solver, the
        standard technique used when the system arises from the discretization
        of an integral equation whose kernel is the Green's function of a PDE is
        to invoke a proxy surface. The proxy surface is used to provide an
        alternative means of representing the action of particular submatrices,
        but its use must be justified and carefully detailed in order to obtain
        the highest accuracy approximation possible. Standard methods usually rely only on heuristics based on Green's identities; here we provide a formal justification and proxy surface compression strategy for a specific integral equation kernel used in exterior scattering.

\end{itemize}
} Details of the above contributions are provided in the subsequent sections of
the manuscript.

\subsection*{Organization}

The paper is organized as follows: Section~\ref{sec:notation} provides a summary
of the notational convention of the paper, and then Section~\ref{sec:setup}
provides a brief derivation of the boundary integral equations of interest and
an overview of the Nystr\"om method for their discretization. In
Section~\ref{sec:skel} we review the notion of strong skeletonization for
matrices, i.e. the algebraic analog of FMM-based
compression~\cite{gr-fmm_3d-1999,Ro-fmm-1985,GrRo-fmm-1987,GrRo-fmm-1988,jiang2022skel}.
This was used by Minden, Ho, Damle and Ying~\cite{minden-str_scel-2017} to
develop the {\em recursive strong skeletonization factorization} method (RS-S),
of which our algorithm is a variant. We then present in Section~\ref{sec:quad}
efficient quadrature coupling and our modifications of strong skeletonization,
leading to the FMM-LU scheme. (Closely related are the inverse FMM
approach~\cite{ambikasaran-2014,coulier-ifmm_prec-2017} and the $\mathcal{H}^2$
matrix formalism~\cite{Borm-h2-2010,Hakb-h2-lib}.) Section~\ref{sec:numerical}
contains several numerical examples demonstrating the asymptotic scaling of our
algorithm, as well as the accuracy obtained in solving realistic boundary value
problems. In Section~\ref{sec:conclusions}, we provide some guidelines regarding
the current generation of fast direct solvers, discuss the limitations of the
present scheme, and speculate about potential avenues for future improvement.

{\note

\section{Notation}
\label{sec:notation}

As with any presentation of hierarchical matrix algorithms, there is
unfortunately a fair amount of notation. We deviate slightly from the notation
used in~\cite{minden-str_scel-2017} due to the extensive discussion around regions
of space which depend on both the octree data structure and the local
quadrature corrections. Below is a summary of the notation used:
\begin{itemize}
  \item $n$ : The overall dimension of the matrix being factorized and inverted;
  total number of nodes in the discretization of an integral equation.
  \item $G$, $K$ : The Green's function for the Helmholtz equation in three
  dimensions and the kernel for the combined field potential operator.
  \item $\cI$, $\cS$, $\cD$, $\cK$ : The identity, single layer, double layer,
  and combined field operators appearing in various integral equation
  formulations, respectively.
  \item $B$, $N$, $F$, $Q$, $P$ : A box in the octree data structure, the boxes
  composing its near field, and the boxes composing its far field, respectively.
  The far field $F$ will further be partitioned as~$F = Q \cup P$, details to
  follow.
  \item $\mtx{A}$ : Matrices will always be denoted using bold sans-serif font.
  \item $\sfI, \sfJ, \sfR, \ldots$ : Normal-weight sans-serif font is used to
  denote a collection of indices, i.e. non-negative integers corresponding to a
  selection of columns, rows, etc.
  \item $\mtx{A}_{\sfI\sfJ}$: The submatrix of $\mtx{A}$ obtained from the rows
  indexed by~$\sfI$ and the columns indexed by~$\sfJ$.
\end{itemize}
In the above notational style, we will, for example, refer to a box in the
octree data structure as~$B$ and associate to it a collection of indices~$\sfB$
corresponding to the indices of points that lie in that box (it will be made
clear when~$\sfB$ refers to the \emph{active} indices in that box, as defined below). Analogously, the
set of indicies corresponding to its near field~$N$ will be denoted by~$\sfN$,
etc. In this situation, serif fonts refer to regions of space and sans-serif
fonts refer to indices of points located in those regions.

}

\section{Problem setup}
\label{sec:setup}

Let~$\Omega$ be a bounded region in~$\bbR^{3}$ with smooth
boundary~$\Gamma = \partial\Omega$. Given a kernel $K$, consider the following
integral equation for~$\sigma$ on the surface~$\Gamma$:
\begin{equation}
\label{eq:inteq0}
\alpha \sigma(\bx) + \int_{\Gamma} K(\bx,\by) \, \sigma(\by) \, da(\by)= f(\bx) \, .
\end{equation}
Here $f$ is a given function on the boundary, $\sigma$ is
an unknown function to be determined, and $\alpha \in \bC$ a
constant.  Such integral equations (and their analogs in the
vector-valued case) naturally arise in the solution of boundary value
problems for the Laplace, Helmholtz, Yukawa, Maxwell, and Stokes
equations, just to name a few. In these settings, the kernel $K$ is typically
related to the Green's function  (or its derivatives) of the
corresponding partial differential equation.

For example, consider the exterior Dirichlet boundary value problem
for the Helmholtz equation with wave number $k$ and boundary data $f$:
the Helmholtz potential $u$ defined in $\mathbb{R}^{3}\setminus
\Omega$ satisfies
\begin{equation}\label{eq:extdir}
  \begin{aligned}
    (\Delta + k^2) \, u &=  0 &\qquad &\text{in } \bbR^{3}
    \setminus \Omega,\\
    u &= f & &\text{on } \Gamma, \\
    \lim_{r \to \infty} \, r \left(\frac{\partial u}{\partial r} - ik u
    \right) &= 0. & &
\end{aligned}
\end{equation}
Let $G = G(\bx,\by)$ denote the free-space Green's function for
the Helmholtz equation given by
\begin{equation}
\label{eq:greenfunhelm}
G(\bx,\by)  = \frac{e^{ik |\bx-\by|}}{4\pi |\bx-\by|} \, ,
\end{equation}
and let 
\begin{equation}
\label{eq:comb_ref}
K(\bx,\by) = (\bn(\by) \cdot \nabla_{\by}G(\bx,\by)) - ik G(\bx,\by) \, ,
\end{equation}
where~$\bn(\by)$ is the outward normal at $\by \in \Gamma$. The above
kernel is the kernel of what is known as the~\emph{combined field
  representation}.
Suppose
that $\sigma$ satisfies~\eqref{eq:inteq0} with $\alpha=1/2$ and $K$ as
defined above, then the potential $u$ given by
\begin{equation}
u(\bx) = \int_{\Gamma} K(\bx,\by) \, \sigma(\by)\, da(\by) \, ,\quad  \bx \in \mathbb{R}^{3}\setminus \Omega \, ,
\end{equation}
is the solution to the exterior Dirichlet problem for the Helmholtz equation
in~\eqref{eq:extdir}. We will often denote the above formula more succinctly
as~$u = \cK_\Gamma \sigma$, where~$\cK_\Gamma$ is the \emph{combined field} layer potential operator
along~$\Gamma$ with kernel~$K$. Note: the operator~$\cK_{\Gamma}$ is the sum of
a \emph{double layer} operator and a \emph{single layer} operator, given by
\begin{equation}
  \begin{aligned}
    \cK_{\Gamma}[\sigma](\bx) &= \cD_{\Gamma}[\sigma](\bx) - ik \cS_{\Gamma}[\sigma](\bx) \\
    &= \int_{\Gamma} (\bn(\by) \cdot \nabla_{\by} G(\bx,\by)) \, \sigma(\by) \, da(\by) - ik \int_{\Gamma} G(\bx,\by) \, \sigma(\by) \, da(\by) .
  \end{aligned}
\end{equation}
More details on layer potential representations of solutions to the Helmholtz
equation and the associated integral equations can be found
in~\cite{kress_2014}, for example.

The integral equation~\cref{eq:inteq0} can be discretized with
high-order accuracy using (for example) a suitable Nystr\"{o}m
method~\cite{atkinson_1997} resulting in the following linear system
\begin{equation}
\label{eq:disc1}
\alpha \sigma_{i} + \sum_{j=1, j\neq i}^{n} K(\bx_{i},\bx_{j}) \, \sigma_{j}
\, w_{ij} = f(\bx_{i}) \, .
\end{equation}
Here $\bx_{i}$ and $w_{ij}$ are the quadrature nodes and weights, respectively,
while $\sigma_{i}$ is an approximation to the true value $\sigma(\bx_{i})$. The
kernel $K=K(\bx,\by)$ is often singular when~$\bx=\by$ and smooth otherwise.
When using high-order discretization methods, it is often possible to use
quadrature weights independent of the target location, except possibly for a
local neighborhood of targets close to each source, i.e. $w_{ij} = w_{j}$ for
all $\bx_{j} \in \textrm{Far}(\bx_{i})$, for a suitably defined well-separated
region $\textrm{Far}(\bx_i)$. Such quadratures are often referred to as
locally-corrected quadratures for obvious reasons. Many existing quadrature
methods such as coordinate transformation methods, singularity subtraction
methods, quadrature by expansion, Erichsen-Sauter rules, and generalized
Gaussian methods combined with adaptive integration can be used as
locally-corrected
quadratures~\cite{bruno2001fast,bruno_garza_2020,malhotra19,erichsen1998quadrature,Siegel2018ALT,Wala2018,Wala2020,ying,bremer_2012c,bremer_2013,bremer-2015,bremer,erichsen1998quadrature,greengard2021fast,wu2020corrected}.
In this setting, the discrete linear system~\cref{eq:disc1} then takes the form
\begin{equation}
\label{eq:disc2}
\alpha \sigma_{i} + \sum_{\substack{j=1\\ \bx_{j} \in
    \textrm{Far}(\bx_{i})}}^{n} K(\bx_{i},\bx_{j}) \, \sigma_{j} \, w_{j}  +
\sum_{\substack{j=1\\ \bx_{j} \not \in \textrm{Far}(\bx_{i})}}^{n}
K(\bx_{i},\bx_{j}) \, \sigma_{j} \, w_{ij}  = f(\bx_{i}) \, .
\end{equation}
We further rescale the above equation so that the unknowns are
instead~$\tilde{\sigma}_{j} = \sigma_{j} \sqrt{w_{j}}$, for which the discrete
linear system becomes
\begin{equation}
\label{eq:disc3}
\alpha \tilde{\sigma}_{i} + \sum_{\substack{j=1\\ \bx_{j} \in
    \textrm{Far}(\bx_{i})}}^{n}  \sqrt{w}_{i} \, K(\bx_{i},\bx_{j})
\, \sqrt{w_{j}} \, \tilde{\sigma}_{j}   + 
\sum_{\substack{j=1\\ \bx_{j} \not \in \textrm{Far}(\bx_{i})}}^{n}
\sqrt{w_{i}} \, K(\bx_{i},\bx_{j}) \, \frac{w_{ij}}{\sqrt{w_{j}}} \,\sigma_{j}   = \sqrt{w_{i}} \, f(\bx_{i}) \, .
\end{equation}
The scaling by the square root of the weights in the above equation
formally embeds the discrete solution~$\sigma_j$ in ${L}^2$
and 
results in the discretized operators (including sub-blocks of the
matrices) to have norms and condition numbers which are close to (and
converge to) those
of the continuous versions of the corresponding
operators~\cite{bremer_2012}.  In the following work, we restrict our
attention to the fast solution of linear systems of the
form~\cref{eq:disc3}.

\begin{remark}
  In many locally-corrected quadrature methods, in order to accurately compute
  far interactions the underlying discretization may require some additional
  oversampling to meet an \emph{a priori} specified precision requirement.
  See~\cite{greengard2021fast} for a thorough discussion. In short, let
  $\bs_{j}$, $\hat{w}_{j}$, $\hat{\sigma}_{j}$, for $j=1,2,\ldots m>n$, denote
  the set of oversampled discretization nodes, the corresponding quadrature
  weights for integrating smooth functions on the surface, and the interpolated
  density at the oversampled discretization nodes, respectively. For
  triangulated/quadrangulated surfaces, the values of the density at the
  oversampled nodes can be obtained via polynomial interpolation of the chart
  information (i.e. parameterization and metric tensor) and the values of the
  density at the original discretization nodes. For example, in the case of a
  quadrilateral patch, these interpolation operations can be computed explicitly
  using tensor-product Legendre polynomials and the matrices which map function
  values to coefficients in a Legendre polynomial expansion.

The
linear system (without the square root scaling) then takes the form
\begin{equation}
\alpha \sigma_{i} + \sum_{\substack{j=1\\ \bs_{j} \in
    \textrm{Far}(\bx_{i})}}^{m} K(\bx_{i},\bs_{j}) \,
\hat{\sigma}_{j} \, \hat{w}_{j}  + 
\sum_{\substack{j=1\\ \bx_{j} \not \in \textrm{Far}(\bx_{i})}}^{n}
K(\bx_{i},\bx_{j}) \, \sigma_{j} \, w_{ij}  = f(\bx_{i}) \, .
\end{equation}
The need for oversampling often arises when using low-order discretizations. The
extension of our approach, and other existing approaches, to discretizations
which require oversampling for far interactions is currently being pursued and
the results will be reported at a later date. The main ideas are similar, but
constructing an efficient algorithm requires carefully addressing many
non-trivial implementation details.
\end{remark}

\section{\note The FMM-LU factorization}
\label{sec:skel}

The basic structure of the FMM-LU factorization is closely related to the
recursive strong skeletonization factorization (RS-S) introduced
in~\cite{minden-str_scel-2017}. In this section, we briefly review key elements
of RS-S. To the extent possible, we use the same notation as
in~\cite{minden-str_scel-2017} to clearly highlight our modifications to their
approach. {\note We provide additional discussion addressing the
  Schur complement update procedure in the case where the algorithm is being
  applied to an adaptive data structure; see Section~\ref{sec:rssadap} below.
  This situation must be considered when computing the far field partitioning of
  points, as will be seen shortly. }

Suppose that all the discretization points $\bx_{i} \in \Gamma$,
$i=1,2,\ldots n$ are contained in a cube $C$.  We will superimpose
on~$C$ a hierarchy of refinements as follows: the root of the tree is
$C$ itself and defined as {\em level 0}.  Level $l+1$ is obtained from
level~$l$ recursively by subdividing each cube at level $l$ into eight
equal parts so long as as the number of points in that cube at level $l$ is
greater than some specified parameter~$s = \cO(1)$.  The eight cubes
created by subdividing a larger cube are referred to as
\emph{children} of the \emph{parent} cube.  When the refinement has
terminated, $C$ is covered by disjoint childless boxes at various
levels of the hierarchy (depending on the local density of the given
points). These childless boxes are referred to as leaf boxes.  For any
box $B$ in the hierarchy, the {\em near field} region of $B$ consists
of other boxes at the same level that touch $B$, and the {\em far
  field} region is the remainder of the domain.
    
For simplicity, we assume that the above described octree satisfies a
standard restriction -- namely, that two leaf nodes which share a
boundary point must be no more than one refinement level apart. In
creating the adaptive data structure as described above, it is very
likely that this level-restriction criterion is not met. Fortunately,
assuming that the tree constructed to this point has $\cO(n)$ leaf
boxes and that its depth is of the order $\cO(\log n)$, it is
straightforward to enforce the level-restriction in a second step
requiring~$\cO(n \log n)$ effort with only a modest amount of
additional refinement~\cite{treebook}.

The near field region and far field region of a box $B$ in the octree hierarchy
is almost always different from the near region and far regions of source and
target locations associated with the locally-corrected quadrature methods
in~\eqref{eq:disc3}. In practice, for almost all targets, the near field region
for the local quadrature corrections is a subset of the near field
region of the leaf box~$B$ containing the target. In the event that this
condition is violated, the RS-S algorithm of~\cite{minden-str_scel-2017} would
require modifications, present in the work of this paper, to handle the
associated quadrature corrections discussed in~\cref{sec:farsplit}.

\subsection{Strong skeletonization}
\label{subsec:strong-skel}

The idea of \emph{strong} skeletonization was recently introduced
in~\cite{minden-str_scel-2017} and extends the idea of using the interpolative
decomposition to globally compress a low-rank matrix to the situation where only
a particular off-diagonal block is low-rank. Effectively, strong skeletonization
decouples some columns/rows of the matrix from the remaining ones. In the
context of solving a discretized boundary integral equation, the off-diagonal
low-rank block is a result of far field interactions via the kernel (e.g.
Green's function) of the integral equation. We briefly recall the standard
interpolative decomposition~\cite{liberty2007randomized} and the form of strong
skeletonization as presented in~\cite{minden-str_scel-2017}.

\begin{definition}[Interpolative decomposition]
Given a matrix $\mtx{A} \in \bC^{|\sfI|\times |\sfJ|}$ with
rows indexed by $\sfI$ and columns indexed by $\sfJ$, an
$\varepsilon$-accurate interpolative decomposition (ID) of $\mtx{A}$ is a
partitioning of $\sfJ$ into a set of so-called skeleton columns denoted
by $\sfS \subset \sfJ$ and redundant columns $\sfR = \sfJ \setminus \sfS$,
and a construction of a corresponding interpolation matrix $\mtx{T}$ such that
\begin{equation}
\left\Vert \mtx{A}_{\sfI \sfR} - \mtx{A}_{\sfI \sfS} \mtx{T}\right \Vert
\leq \varepsilon \Vert \mtx{A}_{\sfI\sfJ} \Vert \, ,
\end{equation}
i.e. the redundant columns are well-approximated, in a relative sense, to the required
tolerance $\varepsilon$ by linear combinations of the skeleton
columns. Equivalently, after an application of a permutation
matrix~$\mtx{P}$ such that $\mtx{A} \mtx{P} = [ \mtx{A}_{\sfI \sfR}
  \quad \mtx{A}_{\sfI \sfS}]$, the interpolative decomposition results
in the $\varepsilon$-accurate low-rank factorization $\mtx{A}_{\sfI \sfJ}
\mtx{P} \approx \mtx{A}_{\sfI \sfS} [\mtx{T} \quad \mtx{I}]$ with the error
estimate,
\begin{equation}
\left\Vert \mtx{A}_{\sfI \sfJ} \mtx{P} - \mtx{A}_{\sfI \sfS} [\mtx{T}
  \quad \mtx{I}] \right\Vert
\leq \varepsilon \Vert \mtx{A}_{\sfI \sfJ} \Vert  .
\end{equation}
The norms above can be taken to be the standard induced spectral norm.
\end{definition}
The ID is most robustly computed using the strong rank-revealing QR
factorization of Gu and Eisenstat~\cite{gu1996efficient}. However, in
this work we use a standard greedy column-pivoted
QR~\cite{martinsson2014}. While both algorithms have similar computational
complexity when $|\sfJ| \leq |\sfI|$, i.e. $\cO( |\sfI| \cdot |\sfJ|^2)$,
the greedy column pivoted QR tends to have better computational
performance.

Next, consider a three-by-three block matrix~$\mtx{A}\in \bC^{n\times n}$ and
suppose that $[n] = \sfI \cup \sfJ \cup \sfK$ is a partition of the index set, with
$[n] = \{ 1, 2, \ldots, n\}$, such that $\mtx{A}_{\sfI \sfK} = \mtx{0}$ and
$\mtx{A}_{\sfK \sfI} = \mtx{0}$, i.e.,
\begin{equation}
\label{eq:block-schur}
  \mtx{A} = \begin{bmatrix}
  \mtx{A}_{\sfI \sfI} & \mtx{A}_{\sfI \sfJ} & \mtx{0}\\
  \mtx{A}_{\sfJ \sfI} & \mtx{A}_{\sfJ \sfJ} &  \mtx{A}_{\sfJ \sfK} \\
 \mtx{0} & \mtx{A}_{\sfK \sfJ} &  \mtx{A}_{\sfK \sfK}
  \end{bmatrix} \, .
\end{equation}
Assuming that $\mtx{A}_{\sfI \sfI}$ is invertible, then using block
Gaussian elimination the matrix $\mtx{A}_{\sfI \sfI}$ can be \emph{decoupled}
from the rest of the matrix as follows
\begin{equation}
\mtx{L} \cdot \mtx{A} \cdot \mtx{U} = 
\begin{bmatrix}
\mtx{I} & \mtx{0} & \mtx{0} \\
-\mtx{A}_{\sfJ \sfI} \mtx{A}_{\sfI \sfI}^{-1} & \mtx{I} & \mtx{0} \\
\mtx{0} & \mtx{0} & \mtx{I}\\
\end{bmatrix} \cdot \mtx{A} \cdot 
\begin{bmatrix}
\mtx{I} & - \mtx{A}_{\sfI \sfI}^{-1} \mtx{A}_{\sfI \sfJ} & \mtx{0} \\
\mtx{0} & \mtx{I} & \mtx{0} \\
\mtx{0} & \mtx{0} & \mtx{I}\\
\end{bmatrix} = 
\begin{bmatrix}
  \mtx{A}_{\sfI \sfI} & \mtx{0} & \mtx{0}\\
\mtx{0} & \mtx{S}_{\sfJ \sfJ} &  \mtx{A}_{\sfJ \sfK} \\
 \mtx{0} & \mtx{A}_{\sfK \sfJ} &  \mtx{A}_{\sfK \sfK}
 \end{bmatrix} \, ,
\end{equation}
where
$\mtx{S}_{\sfJ \sfJ} = \mtx{A}_{\sfJ \sfJ} - \mtx{A}_{\sfJ \sfI}
\mtx{A}_{\sfI \sfI}^{-1} \mtx{A}_{\sfI \sfJ}$ is the only non-zero block
of the matrix that has been modified. The matrix
$- \mtx{A}_{\sfJ \sfI} \mtx{A}_{\sfI \sfI}^{-1} \mtx{A}_{\sfI \sfJ}$ is
often referred to as the Schur complement update.

Suppose now that the matrix $\mtx{A}$ arises from the discretization of an
integral operator. Furthermore, suppose that $\sfB$ is a set of indices of points
$\bx_{i}$ contained in a box $B$ in the octree, that $\sfN$ are the indices of
the set of points contained in the near field region~$N$ of~$B$, and that $\sfF$ are
the indices of points contained in~$B$'s far field region~$F$. Using an
appropriate permutation matrix~$\mtx{P}$, we can obtain the following block
structure for~$\mtx{A}$:
\begin{equation}
\label{eq:proxy-split-0}
\mtx{P}^{T}\mtx{A}\mtx{P} = 
\begin{bmatrix}
\mtx{A}_{\sfB  \sfB} & \mtx{A}_{\sfB \sfN} & \mtx{A}_{\sfB \sfF} \\
\mtx{A}_{\sfN  \sfB} & \mtx{A}_{\sfN \sfN} & \mtx{A}_{\sfN \sfF} \\
\mtx{A}_{\sfF  \sfB} & \mtx{A}_{\sfF \sfN} & \mtx{A}_{\sfF \sfF}
\end{bmatrix} \, .
\end{equation}
The blocks corresponding to interactions between points in $B$ and its far
field, i.e. $\mtx{A}_{\sfB \sfF}$ and $\mtx{A}_{\sfF \sfB}$, are assumed to be
numerically low-rank (which can be seen from an analysis of the underlying PDE or integral equation~\cite{martinsson-book}) and can therefore be compressed using interpolative
decompositions. As before, we partition $\sfB$ into a collection of redundant
points $\sfR$ and a set of skeleton points $\sfS$ such that, up to an appropriate
permutation of rows and columns (which can be absorbed into the permutation
matrix $\mtx{P}$ above), we have
\begin{equation}
\label{eq:proxy-split-1}
\begin{bmatrix}
\mtx{A}_{\sfF \sfB} \\
\mtx{A}_{\sfB \sfF}^{T}
\end{bmatrix} =
\begin{bmatrix}
\mtx{A}_{\sfF \sfR} & \mtx{A}_{\sfF \sfS} \\
\mtx{A}^{T}_{\sfR \sfF} & \mtx{A}^{T}_{\sfR \sfS}
\end{bmatrix}  \approx
\begin{bmatrix}
\mtx{A}_{\sfF \sfS} \\
\mtx{A}^{T}_{\sfS \sfF}
\end{bmatrix} \cdot \begin{bmatrix}
\mtx{T} &  \mtx{I}
\end{bmatrix} \, .
\end{equation}
Note that in the relationship above the same interpolation
matrix~$\mtx{T}$ is used to compress both~$\mtx{A}_{\sfF \sfB}$
and~$\mtx{A}^T_{\sfB \sfF}$. Clearly this is possible when the kernel of the
associate integral equation is symmetric. If the kernel is not
symmetric, empirically the same matrix~$\mtx{T}$ can be used at the
cost of a small increase in~$|\sfS|$. Using the same matrix~$\mtx{T}$
simplifies various subsequent linear algebra manipulations, but
strictly speaking, is not necessary. Different interpolation matrices
can be used for each of~$\mtx{A}^T_{\sfB \sfF}$ and~$\mtx{A}_{\sfF
  \sfB}$. Using different compression matrices will result in smaller
skeleton sets at the cost of more linear algebraic bookkeeping.

Further splitting the indices $\sfB = \sfR \cup \sfS$
in~\cref{eq:proxy-split-0}, and combining
with~\cref{eq:proxy-split-1}, we get
\begin{equation}
  \label{eq:pap}
\mtx{P}^{T}\mtx{A}\mtx{P} = 
\begin{bmatrix}
\mtx{A}_{\sfR  \sfR} & \mtx{A}_{\sfR \sfS} & \mtx{A}_{\sfR \sfN} & \mtx{A}_{\sfR \sfF} \\
\mtx{A}_{\sfS  \sfR} & \mtx{A}_{\sfS \sfS} & \mtx{A}_{\sfS \sfN} & \mtx{A}_{\sfS \sfF} \\
\mtx{A}_{\sfN  \sfR} & \mtx{A}_{\sfN \sfS} & \mtx{A}_{\sfN \sfN} & \mtx{A}_{\sfN \sfF} \\
\mtx{A}_{\sfF  \sfR} & \mtx{A}_{\sfF \sfS} & \mtx{A}_{\sfF \sfN} & \mtx{A}_{\sfF \sfF}
\end{bmatrix}  \approx
\begin{bmatrix}
\mtx{A}_{\sfR  \sfR} & \mtx{A}_{\sfR \sfS} & \mtx{A}_{\sfR \sfN} & \mtx{T}^{T}\mtx{A}_{\sfS \sfF} \\
\mtx{A}_{\sfS  \sfR} & \mtx{A}_{\sfS \sfS} & \mtx{A}_{\sfS \sfN} & \mtx{A}_{\sfS \sfF} \\
\mtx{A}_{\sfN  \sfR} & \mtx{A}_{\sfN \sfS} & \mtx{A}_{\sfN \sfN} & \mtx{A}_{\sfN \sfF} \\
\mtx{A}_{\sfF \sfS} \mtx{T} & \mtx{A}_{\sfF \sfS} & \mtx{A}_{\sfF \sfN} & \mtx{A}_{\sfF \sfF}
\end{bmatrix} \, .
\end{equation}
Since the redundant rows and columns of the interaction between points in
$\sfR$ and $\sfF$ can be well-approximated by the
corresponding rows and columns of the interactions between points in $\sfS$
and $\sfF$, we can decouple points $\sfR$ from the far field
points $\sfF$ as follows. Let $\mtx{E}$, $\mtx{F}$ denote the
elimination matrices defined on the partition $[n] = \sfR \cup \sfS \cup
\sfN \cup \sfF$ as
\begin{equation}
\mtx{E} = \begin{bmatrix}
\mtx{I} & -\mtx{T}^{T} & & \\
& \mtx{I} & & \\
& & \mtx{I} & \\
& & & \mtx{I}\end{bmatrix} \qquad
\text{and} \qquad
\mtx{F} = \begin{bmatrix}
\mtx{I} & & & \\
-\mtx{T} & \mtx{I} & & \\
& & \mtx{I} & \\
& & & \mtx{I}\end{bmatrix} \, .
\end{equation} 
Then
\begin{equation}
\mtx{E} \mtx{P}^{T} \mtx{A} \mtx{P} \mtx{F} =
\begin{bmatrix}
\mtx{X}_{\sfR  \sfR} & \mtx{X}_{\sfR \sfS} & \mtx{X}_{\sfR \sfN} & \mtx{0} \\
\mtx{X}_{\sfS  \sfR} & \mtx{A}_{\sfS \sfS} & \mtx{A}_{\sfS \sfN} & \mtx{A}_{\sfS \sfF} \\
\mtx{X}_{\sfN  \sfR} & \mtx{A}_{\sfN \sfS} & \mtx{A}_{\sfN \sfN} & \mtx{A}_{\sfN \sfF} \\
\mtx{0} & \mtx{A}_{\sfF \sfS} & \mtx{A}_{\sfF \sfN} & \mtx{A}_{\sfF \sfF}
\end{bmatrix} \, ,
\end{equation}
where the notation $\mtx{X}_{\sfI \sfJ}$ is used to indicate blocks of
the above matrix whose entries are different from the entries of the
original matrix in~\eqref{eq:pap}.  We also note that the matrices
$\mtx{E}$ and $\mtx{F}$ are, in fact, block diagonal when viewed over
the partition $[n] = \sfB \cup \sfN \cup \sfF$, and therefore the above
factorization can be considered a type of $\mtx{LU}$ elimination.

Now, assuming again that the block $\mtx{X}_{\sfR \sfR}$ is invertible, we can use
it as a pivot block to completely decouple the redundant indices $\sfR$ from
the rest of the problem as follows. Let $\mtx{L}$ and $\mtx{U}$ denote
the block upper and lower triangular elimination matrices given by
\begin{equation}
\mtx{L} = 
\begin{bmatrix}
\mtx{I} & & & \\
-\mtx{X}_{\sfS \sfR} \mtx{X}_{\sfR \sfR}^{-1} & \mtx{I} & & \\
-\mtx{X}_{\sfN \sfR} \mtx{X}_{\sfR \sfR}^{-1} & & \mtx{I} & \\
& & & \mtx{I} \\
\end{bmatrix} \qquad \text{and} \qquad
\mtx{U} = 
\begin{bmatrix}
\mtx{I} & -\mtx{X}_{\sfR \sfR}^{-1} \mtx{X}_{\sfR \sfS} & -\mtx{X}_{\sfR \sfR}^{-1} \mtx{X}_{\sfR \sfN} & \\
& \mtx{I} & & \\
& & \mtx{I} & \\
& & & \mtx{I} \\
\end{bmatrix} \, .
\end{equation}
Then, we set
\begin{equation}
  \label{eq:strong-skel-elim}
  \begin{aligned}
\mtx{Z}\left(\mtx{A}; \sfB \right)  & = \mtx{L} \mtx{E} \mtx{P}^{T}
\mtx{A} \mtx{P} \mtx{F} \mtx{U} \\
&= 
\begin{bmatrix}
\mtx{X}_{\sfR  \sfR} & \mtx{0} & \mtx{0} & \mtx{0} \\
\mtx{0} & \mtx{X}_{\sfS \sfS} & \mtx{X}_{\sfS \sfN} & \mtx{A}_{\sfS \sfF} \\
\mtx{0} & \mtx{X}_{\sfN \sfS} & \mtx{X}_{\sfN \sfN} & \mtx{A}_{\sfN \sfF} \\
\mtx{0} & \mtx{A}_{\sfF \sfS} & \mtx{A}_{\sfF \sfN} & \mtx{A}_{\sfF \sfF}
\end{bmatrix}.
\end{aligned}
\end{equation}
The matrix in~\cref{eq:strong-skel-elim} is of the
form~\cref{eq:block-schur}.  This process of decoupling the redundant
degrees of freedom in~$\sfB$ from the rest of the problem is referred to
as the \emph{strong skeletonization of~$\mtx{A}$ with respect
  to~$\sfB$}.
The resulting matrix is denoted by~$\mtx{Z}\left(\mtx{A};
\sfB \right)$,
as above.

Since the matrices $\mtx{E}$ and $\mtx{F}$ are block diagonal with
respect to the partition $[n] = \sfB \cup \sfN\cup \sfF$,
equation~\eqref{eq:strong-skel-elim} can be re-written in terms of a
block $\mtx{LU}$-like factorization of the original matrix~$\mtx{A}$:
\begin{equation}
    \mtx{A} = \left( \mtx{P} \mtx{E}^{-1} \mtx{L}^{-1} \right) \mtx{Z}(\mtx{A}; \sfB) \left( \mtx{U}^{-1} \mtx{F}^{-1}\mtx{P}^{T} \right) \, .
\end{equation}
For notational convenience, let $\mtx{V}$ and $\mtx{W}$ denote the
left and right skeletonization operators defined by
\begin{equation}
\label{eq:lrskeldef}
\mtx{V} = \mtx{P} \mtx{E}^{-1} \mtx{L}^{-1}, \qquad  \mtx{W}
= \mtx{U}^{-1} \mtx{F}^{-1} \mtx{P}^{T}  ,
\end{equation}
with the understanding that these matrices will be stored and used in
factored form for computational efficiency. Moreover, the matrices
$\mtx{L},\mtx{U}, \mtx{E}^{-1}, \mtx{F}^{-1}$ are block triangular
matrices with identities on the diagonal and hence their inverses can
be trivially computed by toggling the sign of the nonzero off-diagonal
blocks.  With this shorthand, we can obtain an even
more compact representation of~$\mtx{Z}(\mtx{A}; \sfB)$ given by
\begin{equation}
\mtx{Z}(\mtx{A}; \sfB) = \mtx{V}^{-1} \mtx{A} \mtx{W}^{-1} \, .
\end{equation}

\begin{remark}
  The elimination matrices $\mtx{E}$ and $\mtx{F}$ are referred to as
  $\mtx{U}_{T}$ and $\mtx{L}_{T}$ in~\cite{minden-str_scel-2017}. However, for
  clarity of identifying the block lower and upper triangular structure
  in~\mbox{$\mtx{V} = \mtx{P} \mtx{E}^{-1} \mtx{L}^{-1}$} and~\mbox{$\mtx{W}=\mtx{U}^{-1} \mtx{F}^{-1} \mtx{P}^{T}$} with respect to the partition $[n]
  = \sfB \cup \sfN \cup \sfF$, we have renamed $\mtx{U}_{T} =
  \mtx{E}$ and $\mtx{L}_{T} = \mtx{F}$.
\end{remark}


\subsection{Recursive strong skeletonization}
\label{subsec:rs-s-alg}

In this section, we provide a short summary of the RS-S algorithm. We refer the
reader to the original manuscript for a more detailed
description~\cite{minden-str_scel-2017}. The RS-S algorithm proceeds by
sequentially applying the strong skeletonization procedure to each box in the
level-restricted tree, where, as mentioned before it is assumed that that two
leaf nodes which share a boundary point must be no more than one refinement
level apart. The boxes in the tree hierarchy are traversed in an upward pass,
i.e. boxes at the finest level will be processed first followed by boxes at
subsequent coarser levels. After each application of the strong skeletonization
procedure, only the skeleton points $\sfS$ associated with each box are retained
for further processing. These will be referred to as the \emph{active degrees of
  freedom}. Even when constructing near field and far field index sets of a box,
only the active degrees of freedom contained in the respective regions are
retained (other degrees of freedom have been deemed \emph{redundant} and
decoupled from the system). For boxes at coarser levels, the active degrees of
freedom for each box is the union of the active degrees of freedom of each of
its children boxes. After regrouping the active indices from all the children of
boxes at coarser levels, the process of strong skeletonization can be applied to
those boxes as well. This process is continued until there are no remaining
active degrees of freedom in the far field region of any box at a given level or
the algorithm reaches level 1 in the tree structure (for which the statement is
trivially true since there are no boxes in the far field region of any box).

 Similar to before, let $\mtx{V}_{i}$ and $\mtx{W}_{i}$ denote the
 left and right skeltonization operators associated with box $B_{i}$
 defined in~\cref{eq:lrskeldef}. Suppose that the multi-level RS-S
 algorithm of~\cite{minden-str_scel-2017} terminates at box~$B_{M}$;
 let~$\sfB_{t}$ denote the remaining active degrees of freedom in the
 domain. Let $\mtx{P}_{t}$ denote the permutation which orders the
 points in $\sfB_{t}$ in a contiguous manner. Then the RS-S
 factorization of the matrix $\mtx{A}$ takes the form
\begin{equation}
\mtx{A} \approx \left( \mtx{V}_{1} \mtx{V}_{2} \cdots \mtx{V}_{M} \right) \mtx{P}_{t} \mtx{D} \mtx{P}_{t}^{T} \left( \mtx{W}_{M} \mtx{W}_{M-1} \cdots \mtx{W}_{1}  \right) \, .
\end{equation}
Here $\mtx{D}$ is the block diagonal matrix given by
\begin{equation}
\mtx{D} = \begin{bmatrix}
\mtx{X}_{\sfR_{1} \sfR_{1}} & &  \\
& \ddots & &  \\
& & \mtx{X}_{\sfR_{M} \sfR_{M}} & \\
& & & \mtx{A}_{\sfB_{t} \sfB_{t}}
\end{bmatrix} \, ,
\end{equation}
where $\sfR_{j}$ are the redundant indices in box~$B_{j}$. An approximate
factorization of $\mtx{A}^{-1}$ is readily obtained from the formula above and
is given by
\begin{equation}
\mtx{A}^{-1} \approx \left(\mtx{W}_{1}^{-1} \cdots \mtx{W}_{M}^{-1} \right) \mtx{P}_{t} \mtx{D}^{-1} \mtx{P}_{t}^{T} \left(\mtx{V}_{M}^{-1} \cdots \mtx{V}_{1}^{-1} \right) \, .
\end{equation}
In the event that the matrix $\mtx{A}$ is positive definite, one can also
compute the generalized square-root and the
log-determinant~\cite{minden2017gp,ambikasaran_2016,ambikasaran2016symmetric} of
the matrix~$\mtx{A}$ using the factorization above.

Under mild assumptions on the ranks of the interactions between a box and its
far field, the cost of applying the compressed operator and its inverse
$(t_{a},t_{s})$, and the memory required ($m_{f}$) in the algorithm all scale
linearly in the number of discretization points $n$ independent of the ambient
dimension of the data. {\note Naively computing the low-rank strong
  skeletonizations would incur an overall cost in computing the RS-S
  factorization $(t_{f})$ proportional to~$\cO(n^{2})$. This can be reduced
  to~$\cO(n)$ by compressing via \emph{proxy surfaces} for many important
  classes of interaction kernels, as described in the next section. } In
particular, if the ranks of the far field blocks are $\cO(L-\ell+c)^q$ for boxes
on level~$\ell$ where $c,q>0$, then $t_{f}, t_{s} ,t_{a}, m_{f} = \cO(n)$ with
the implicit constant being a polynomial power of $\log(1/\varepsilon)$, where
$\varepsilon$ is the requested tolerance. These conditions are typically
satisfied for discretizations of boundary integral equations for non-oscillatory
problems.

\subsection{Low-rank approximation using proxy surfaces}
\label{subsec:proxy}

{\note For any box $B$, there are typically~$\cO(1)$ associated active degrees
of freedom (as described in the previous section) and typically $\cO(n)$ points
in its far field region~$F$ (particularly on the finest level of the octree).}
Due to this fact, any algorithm which requires dense assembly of all
blocks~$\mtx{A}_{\sfF \sfB}$ and~$\mtx{A}_{\sfB \sfF}$ for all boxes~$B$ in
order to compute the interpolative decompositions in~\cref{eq:proxy-split-1}
would result in an~$\cO(n^2)$ computational complexity for constructing a
compressed representation of~$\mtx{A}$. In this section, we briefly discuss an
indirect, efficient, and provably accurate approach for constructing the
interpolative decomposition of~$\mtx{A}_{\sfF \sfB}$ and~$\mtx{A}_{\sfB \sfF}$.
More specific implementation details on this accelerated compression in the case
of the combined field kernel are given in Section~\ref{sec:quad}, and general
overviews of the approach can be found
in~\cite{martinsson-book,xing2020proxy,ye2020proxy}.

In order to achieve the linear-time speedup in compression of the off-diagonal
blocks~$\mtx{A}_{\sfF \sfB}$ and~$\mtx{A}_{\sfB \sfF}$, {\note and subsequent
recursive compression with respect to active degrees of freedom at higher levels
in the octree} we will make two additional assumptions regarding the discrete
linear system~\cref{eq:disc3} {\note which are slightly more restrictive that
the assumptions made in~\cite{minden-str_scel-2017}, but which apply directly to
our scattering problem at hand. Alternative assumptions are required in the case
of proxy compression for different boundary value problems}:
\begin{enumerate}
  \item For each box~$B$, there exists a partition of the far field index set $F
        = Q \cup P$, with~\mbox{$|\sfQ| = \cO(1)$}, such that for all ${j} \in
        \sfB$ and ${i} \in \sfP$, the corresponding matrix entries in
        $\mtx{A}_{\sfP \sfB}$ are given by $\sqrt{w}_{i}K(\bx_{i},\bx_{j})
        \sqrt{w}_{j}$. We also assume that the converse holds: for all ${j} \in
        \sfP$ and ${i} \in \sfB$, the corresponding matrix entries in
        $\mtx{A}_{\sfB \sfP}$ are given by $\sqrt{w}_{i} K(\bx_{i},\bx_{j})
        \sqrt{w}_{j}$. This is equivalent to assuming that for these sets of
        well-separated points in~$B$ and~$P$, a smooth source-dependent
        quadrature rule can be used to accurately discretize the underlying
        integral equation. One could, in principle, require this assumption to
        hold for all points in $F$. However, due to Schur complement updates
        arising from the recursive strong skeletonization procedure, as we will
        see, this condition will be violated for some points in~$F$. This
        issue is discussed in further detail in~\cref{subsec:rs-s-alg}.

  \item {\note The kernel~$K = K(\by,\bx)$ is a linear combination of the
        Green's function of a homogeneous elliptic partial differential equation
        (denoted by $G$) and its derivatives in the~$\bx$-variable only
        (extensions to other classes of kernels will be reported at a later
        date). For example, when compressing~$\mtx{A}_{\sfF \sfB}$, if~$\gamma$ is
        a smooth surface embedded in~$\bbR^3$ denoting the proxy surface which
        encloses box~$B$, then the interaction
        kernel~$K$ satisfies the  following conditions, both of
        which will be taken advantage of:
  \begin{itemize}
    \item For $\by \in \gamma$, integrals of the kernel $K(\by,\bx)$ in~$\bx$
          along~$\gamma$ should be interpreted in a principal value sense.
    \item $K(\cdot,\bx)$ satisfies the underlying PDE (e.g. the homogeneous
          Helmholtz equation in this case) at all points except~$\bx$.
  \end{itemize}
        These assumptions are trivially true for both the Helmholtz single and
        double layer potentials, for example.

      }
\end{enumerate}
Alternative assumptions or compression techniques must be made when using
integral representations other than the combined field formulation used in this
work (details on these types of proxy compressions will be reported at a later
date in a subsequent manuscript). The above assumptions will be used in
constructing efficient compression using Green's identities and what are known
as proxy surfaces.

We now give a brief overview of the proxy compression procedure for the
block~$\mtx{A}_{\sfF \sfB}$ which corresponds to what is known as choosing
\emph{outgoing skeletons}. This block of the discretized integral equation maps
source charges and dipoles in~$B$ to potentials in~$F$. In particular,
sources in~$B$ induce a potential -- outside of~$B$ -- which satisfies the
underlying homogeneous elliptic PDE, and therefore, this potential can be
represented by an equivalent charge density distributed along a \emph{proxy
surface}~$\gamma$ which encloses~$B$ but does not include points 
in~$P$. This effectively means that the block~$\mtx{A}_{\sfF
\sfB}$ can be split and factored as
\begin{equation}
  \label{eq:prox1}
  \begin{aligned}
  \mtx{A}_{\sfF \sfB} &\approx
  \begin{bmatrix}
    \mtx{A}_{\sfQ \sfB} \\
    \sqrt{\mtx{D}_\sfP} \mtx{M}_{\sfP \gamma} \mtx{K}_{\gamma \sfB} \sqrt{\mtx{D}_\sfB}
  \end{bmatrix} \\
  &=\begin{bmatrix}
  \mtx{I} & \mtx{0} \\
  \mtx{0} & \sqrt{\mtx{D}_\sfP}\mtx{M}_{\sfP \gamma}
  \end{bmatrix}
  \begin{bmatrix}
     \mtx{A}_{\sfQ \sfB} \\
     \mtx{K}_{\gamma \sfB} \sqrt{\mtx{D}_\sfB}
  \end{bmatrix},
  \end{aligned}
\end{equation}
where~$\mtx{D}_{\sfP}$ and~$\mtx{D}_{\sfB}$ are diagonal matrices which contain
smooth quadrature weights for points in~$P$ and~$B$,
respectively,~$\mtx{K}_{\gamma \sfB}$ is a matrix with entries generated by the
kernel~$K$,
\[
  \mtx{K}_{\gamma \sfB}(i,j) = K(\by_{i},\bx_{j}),
\]
where~$\by_{i}$ are points on the proxy surface~$\gamma$ and~$\bx_{j} \in B$,
and lastly,~$\mtx{M}_{\sfP \gamma}$ is a matrix which maps discrete densities
on~$\gamma$ to potentials in~$F$. Specific details of computing the above
factorization (and of the form of~$\mtx{M}_{\sfP \gamma}$) are contained in
Section~\ref{sec:proxy-fix}, as well as a justification for the existence  of
such a factorization. {\note In~\cite{minden-str_scel-2017} only the symmetric
kernel case of~$\log|x-y|$ was discussed, and the case of Helmholtz and issues
arising from internal resonances was not mentioned.} The implication of the
above factorization is the following: since the number of active degrees of
freedom in~$B$ and~$|\sfQ|$ are~$\cO(1)$, and as will be shown later on~$\gamma$
can be discretized using a modest number of points depending on the diameter of
the box~$B$ in the oscillatory regime, the matrix on the right on the second
line in~\eqref{eq:prox1} has dimensions which are~$\cO(1)$. This means that an
interpolative decomposition on the columns of this matrix can be computed
with~$\cO(1)$ work. Denoting this compression as
\begin{equation}
  \begin{bmatrix}
     \mtx{A}_{\sfQ \sfB} \\
     \mtx{K}_{\gamma \sfB} \sqrt{\mtx{D}_\sfB}
  \end{bmatrix} \mtx{P} \approx
  \begin{bmatrix}
     \mtx{A}_{\sfQ \sfS} \\
     \mtx{K}_{\gamma \sfS} \sqrt{\mtx{D}_{\sfS}}
  \end{bmatrix}
  \begin{bmatrix}
    \mtx{T}_{\sfS \sfR} & \mtx{I}
  \end{bmatrix},
\end{equation}
we have effectively computed a low-rank approximation of~$\mtx{A}_{\sfF \sfB}$ at
a cost of only~$\cO(1)$ flops. The matrix~$\mtx{P}$ above is a permutation
matrix that appropriately re-orders the columns according to those which have
been chosen as \emph{skeleton} columns and those which are \emph{redundant}
columns; the skeleton columns are denoted by the set~$\sfS$ and the redundant
columns are denoted by~$\sfR$. This factorization is worked out in detail later
on and results in equation equation~\eqref{eq:incsr}.

The compression of the dual matrix~$\mtx{A}_{\sfB \sfF}$ can be done in a nearly
identical manner after observing that all potentials which satisfy the
underlying PDE in the box~$B$ can be represented using a charge density lying on
the same proxy surface~$\gamma$, and therefore an interpolative decomposition on
the rows of this matrix can be computed in~$\cO(1)$ time. More details, and a
justification, are found in Section~\ref{sec:proxy-fix}. Lastly, these two
compressions can be done simultaneously in order to generate the same
interpolation~``$\mtx{T}$'' matrix for both~$\mtx{A}_{\sfF \sfB}$
and~$\mtx{A}_{\sfB \sfF}$ at a modest increase in rank. While not strictly
necessary, this somewhat simplifies the subsequent factorization procedure and
associated storage costs.

\begin{remark}
  In~\cite{ho-2014}, the authors discuss compression via the proxy method
  applied to integral equations of the form
\begin{equation}
\label{eq:inteq-comp}
\alpha(\bx) \sigma(\bx) + b(\bx) \int_{\Gamma} K(\bx - \by) c(\by) \sigma(\by) \, da(\by)= f(\bx) \, .
\end{equation}
In particular, when constructing $\mtx{K}^{\ext}$, they emphasize the
need for including both matrices~$\mtx{G}_{\sfP \sfB} \mtx{b}_{\sfB}$ and~$\mtx{G}_{\Gamma
  \sfB} \mtx{c}_{\sfB}$, where $\mtx{b}_{\sfB}$ and $\mtx{c}_{\sfB}$ are diagonal matrices
with entries $b(\bx_{i})$ and $c(\bx_{i})$, for $\bx_{i} \in B$,
respectively. Handling non-uniform quadrature weights $w_{j}$ is
equivalent to compressing a discretized version
of~\cref{eq:inteq-comp} with $b(\bx_{j}) = \sqrt{w_{j}}$, $\alpha(\bx)
= \alpha$, and $c(\bx_{j}) = \sqrt{w_{j}}$.
\end{remark}

\subsection{Adaptive data structures}
\label{sec:rssadap}

One crucial detail still remains to be resolved for the construction of this
factorization: how best to obtain a partition $\sfF = \sfQ \cup \sfP$ of the
active far field indices for each box such that the first condition
in~\cref{subsec:proxy} is satisfied. In particular, it must be shown that even
after several recursive applications of the strong skeletonization procedure the
matrix entries corresponding to $\mtx{A}_{\sfP \sfB}$ are still given by
$\sqrt{w_i} K(\bx_{i},\bx_{j}) \sqrt{w_{j}}$ and are unaffected by Schur
complements introduced during the elimination procedure. {\note To put this another
way: on every level of the recursive algorithm, the entries of~$\mtx{A}_{\sfP \sfB}$ must be pure kernel evaluations and
unaffected by any local quadrature corrections or Schur complement updates so
that compression can be performed via a proxy surface (as discussed in the
previous subsection).}  That this is even
possible is not immediately obvious since the Schur complement update obtained
by applying the strong skeletonization procedure to one box might potentially
affect the interaction between a different box and its far interactions. Recall
that the Schur complement constructed during the strong skeletonization
procedure applied to box $B$ updates the interaction between points in the near
field region of~$B$. Referring to Figure 5 in~\cite{minden-str_scel-2017},
reproduced here in Figure~\ref{fig:minden}, using strong skeletonization to
compress the operator $\mtx{A}$ with respect to points in $B_{1}$ updates the
entries of the matrix corresponding to interactions of points contained in
$B_{2}$ and $B_{3}$ since both of them are in the near field region of $B_{1}$.
However, $B_{3}$ is in the far field region of $B_{2}$, and thus the entries
of~$\mtx{A}_{\sfF_{2} \sfB_{2}}$ (where $F_{2}$ denotes the far field region of
$B_{2}$) will not correspond to the original matrix entries. However, all such
interactions can be included in the partition~$\sfQ_{2}$ of $\sfF_{2} = \sfQ_{2}
\cup \sfP_{2}$. A graphical depiction of the partitioning in the case where the
tree is an adaptive one is shown in Figure~\ref{fig:adaptive}.

\afterpage{

\begin{figure}[h]
  \centering
  \begin{subfigure}[t]{.45\linewidth}
    \centering
    \begin{tikzpicture}[scale=.7]
      \draw[fill=orange,very thick] (2.0,1.0) rectangle (7.0,6.0);
      \draw[fill=teal,very thick] (3.0,2.0) rectangle (6.0,5.0);
      \draw[fill=lime,very thick] (4.0,3.0) rectangle (5.0,4.0);
      \draw[step=1.0,black,thin] (-0.5,-1.5) grid (8.5,6.5);
      \draw[black,thick] (-0.5,-1.5) rectangle (8.5,6.5);
      \node[anchor=center] at (4.5,3.5) {$B_{1}$};
      \node[anchor=center] at (3.5,2.5) {$B_{2}$};
      \node[anchor=center] at (3.5,4.5) {$B_{3}$};
      \draw[black,very thick,densely dashed] (4.5,3.5) circle (2.5);
    \end{tikzpicture}
    \caption{An initial skeletonization with respect to~$B_1$ updates
      interactions with~$B_2$ and~$B_3$.}
  \end{subfigure} \quad
  \begin{subfigure}[t]{.45\linewidth}
    \centering
    \begin{tikzpicture}[scale=.7]
      \draw[fill=orange,very thick] (1.0,0.0) rectangle (6.0,5.0);
      \draw[fill=teal,very thick] (2.0,1.0) rectangle (5.0,4.0);
      \draw[fill=lime,very thick] (3.0,2.0) rectangle (4.0,3.0);
      \draw[step=1.0,black,thin] (-0.5,-1.5) grid (8.5,6.5);
      \draw[black,thick] (-0.5,-1.5) rectangle (8.5,6.5);
      \node[anchor=center] at (4.5,3.5) {$S_{1}$};
      \node[anchor=center] at (3.5,2.5) {$B_{2}$};
      \node[anchor=center] at (3.5,4.5) {$B_{3}$};
      \draw[black,very thick,densely dashed] (3.5,2.5) circle (2.5);
      \end{tikzpicture}
      \caption{Next, the previously computed Schur complement~$S_1$
        corresponding to the skeletonization with respect to~$B_1$ must be
        included when skeletonizing with respect to~$B_2$ since it is not in its
        far field.}
  \end{subfigure}
  \caption{Successive sequential partitioning of the near and far fields for
    uniform boxes. Teal denotes the near field~$N$
      and orange denotes~$Q \subset F$.}
  \label{fig:minden}
\end{figure}
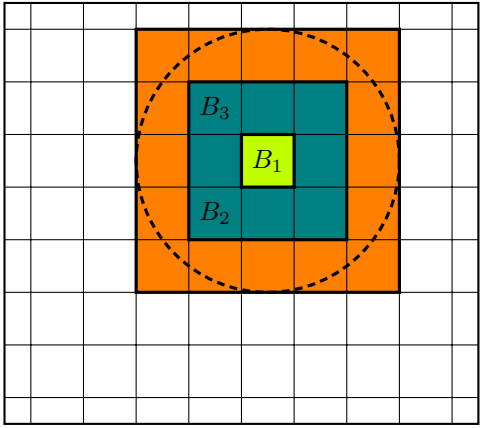
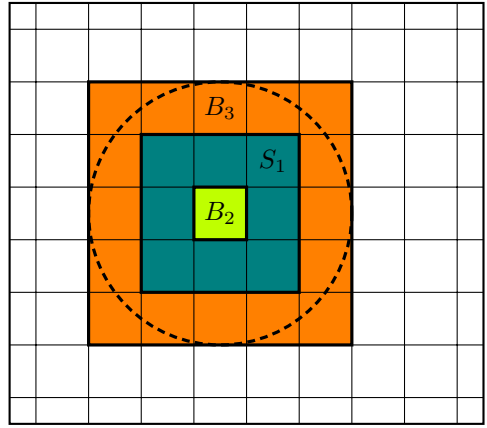

\begin{figure}[h]
  \centering
  \begin{subfigure}[t]{.45\linewidth}
    \centering
    \begin{tikzpicture}[scale=.7]
      \clip (-0.55,-1.55) rectangle (8.55,6.55);
      \draw[step=1.0,black,thin] (-0.5,-1.5) grid (8.5,6.5);
      \draw[fill=orange,very thick] (0.0,-1.5) rectangle (8.5,6.5);
      \draw[step=1.0,black,thin] (0.0,-1.5) grid (8.5,6.5);
      \draw[fill=teal,very thick] (2.0,0.0) rectangle (8.0,6.0);
      \draw[step=1.0,black,thin] (2.0,0.0) grid (8.0,6.0);
      \draw[fill=lime,very thick] (4.0,2.0) rectangle (6.0,4.0);
      \draw[black,thick] (-0.5,-1.5) rectangle (8.5,6.5);
      \node[anchor=center] at (5,3) {$B_{1}$};
      \node[anchor=center] at (3.5,2.5) {$B_{2}$};
      \node[anchor=center] at (3.5,4.5) {$B_{3}$};
      \draw[black,very thick,densely dashed] (5,3) circle (5);
    \end{tikzpicture}
    \caption{Skeletonization on an adaptive tree.}
  \end{subfigure} \quad
  \begin{subfigure}[t]{.45\linewidth}
    \centering
    \begin{tikzpicture}[scale=.7]
      \draw[step=1.0,black,thin] (-0.5,-1.5) grid (8.5,6.5);
      \draw[fill=orange,very thick] (1.0,0.0) rectangle (6.0,5.0);
      \draw[step=1.0,black,thin] (1.0,0.0) grid (6.0,5.0);
      \draw[fill=teal,very thick] (2.0,1.0) rectangle (5.0,4.0);
      \draw[step=1.0,black,thin] (2.0,1.0) grid (5.0,4.0);
      \draw[fill=teal,thin] (4.0,2.0) rectangle (6.0,4.0);
      \draw[very thick] (4.0,4.0) -- (6.0,4.0) -- (6.0,2.0) -- (5.0,2.0);
      \draw[fill=lime,very thick] (3.0,2.0) rectangle (4.0,3.0);
      \draw[black,thick] (-0.5,-1.5) rectangle (8.5,6.5);
      \node[anchor=center] at (5,3) {$S_{1}$};
      \node[anchor=center] at (3.5,2.5) {$B_{2}$};
      \node[anchor=center] at (3.5,4.5) {$B_{3}$};
      \draw[black,very thick,densely dashed] (3.5,2.5) circle (2.5);
      \end{tikzpicture}
      \caption{An adjusted near and far field partition based on the adaptivity.}
  \end{subfigure}
  \caption{An example of successive sequential partitioning of the near and far
    fields with some level-restricted adaptivity. As before, teal denotes the near field~$N$
      and orange denotes~$Q \subset F$.}
  \label{fig:adaptive}
\end{figure}
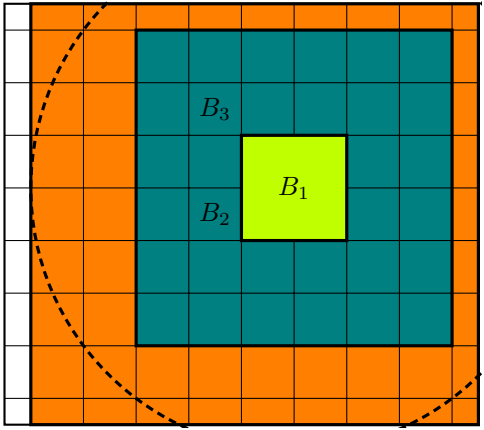
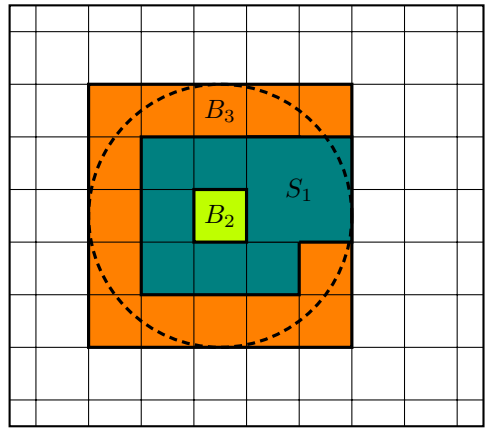

  \clearpage
}

A systematic way of addressing this issue and constructing the partition $\sfF =
\sfQ \cup \sfP$ was presented in~\cite{minden-str_scel-2017}. Suppose that $D$
is a sphere enclosing the box $B$ with radius equal to $5 R/2$, where $R$ is the
side-length of box $B$. Let $\sfP$ denote the set of indices in $\sfF$ such that
if $\bx_{i} \in P \cap B_{\ell}$ (i.e. $i \in \sfP \cap \sfB_\ell$) then $B_{\ell} \in D^c$.
This choice ensures that all the matrix entries in $\mtx{A}_{\sfP \sfB}$ and its
transpose are always of the form $\sqrt{w_{i}} K(\bx_{i},\bx_{j}) \sqrt{w_{j}}$
and $\sqrt{w_{i}} K(\bx_{j},\bx_{i}) \sqrt{w_{j}}$, respectively. Furthermore,
this choice also ensures that after merging the active degrees of freedom from
boxes at finer levels, the partition $\sfF = \sfQ \cup \sfP$ still satisfies the
constraint. These result are proven in~\cite{minden-str_scel-2017} (Theorem 3.1,
and Corollary 3.2 respectively). The additional buffer of placing the proxy
surface separated by two boxes at the same level is crucial in proving those
results.


\begin{remark}
  {\note As detailed in Section 3, in particular Section 3.3.4,
  of~\cite{minden-str_scel-2017}, under certain assumptions on the distribution
  of points and the associated octree data structure, the inversion of the
  system matrix using recursive strong skeletonization scales
  as~$\mathcal O(n)$. In the case of a uniform octree (i.e. all leaf boxes are
  of the same size and on the same level in the hierarchy), it is relatively
  straightforward to track submatrices which have been updated via Schur
  complements. However, if the octree is an adaptive one, and furthermore if it
  contains \emph{quadrature corrections} which are not local to leaves only,
  i.e. that spill out across neighboring blocks which may or may not be ordered
  adjacently in the tree, some additional bookkeeping is needed so as to avoid
  needlessly scanning all possible Schur complement updates.}

  A naive implementation of looping over all Schur complements for updating
  matrix entries could result in an~$\cO(n^2)$ complexity for constructing the
  RS-S factorization (since on the finest level in the data structure, there are
  $\mathcal O(n)$ boxes and therefore $\mathcal O(n)$ associated Schur complements).
  Letting~$\mtx{S}^{(i)}$ denote the Schur complement obtained from eliminating
  the redundant degrees of freedom in box $B_{i}$, one can precompute
  pairs~$(i,j)$ corresponding to Schur complements~$\mtx{S}^{(i)}$ which would
  impact the matrix entries of~$\mtx{A}_{\sfB_{j} (\sfQ_{j} \cup \sfN_{j})}$ and
  its transpose. The need for this is directly due to the additional splitting
  of the far field into~$\sfQ$ and~$\sfP$. Using this list, one can avoid having
  to loop over all Schur complement blocks and therefore retain the~$\cO(n)$
  complexity for constructing the RS-S factorization.
\end{remark}

\section{Quadrature coupling in multiscale
  geometries}
\label{sec:quad}

The primary source of error in constructing the RS-S factorization is
in the compression of the matrix $\mtx{A}_{\sfF \sfB}$  corresponding to far
interactions (and its dual) using the proxy
method~\cite{ye2020proxy,xing2020proxy}.  Recall that after using the proxy method for compression as
discussed in~\cref{subsec:proxy}, the interpolative decomposition
is computed for a matrix whose entries are kernel
evaluations~$K(\bx_i, \bx_j)$. In the following
section, we present three modifications to the standard RS-S algorithm
of~\cite{minden-str_scel-2017}:
\begin{enumerate}
\item Properly formulating the proxy surface compression procedure
  when the kernel~$K$ is obtained using the \emph{combined field
  integral representation};
\item Determining a properly sampled discretization of the proxy
  surface~$\gamma$ using~$n_{\gamma}$ points which is based on the box
size in the tree hierarchy, so as to sufficiently sample the operator when the kernel~$K$ is oscillatory without excessive oversampling; and
\item Finally, constructing a partition $F = Q \cup P$ capable of
handling near field quadrature corrections for multiscale geometries.
\end{enumerate}
We first turn to the accelerated proxy
compression of matrix blocks obtained when using
a combined field representation for the solution to a Dirichlet problem.

\subsection{Proxy compression for combined field representations}
\label{sec:proxy-fix}

As discussed in the introduction, in order to solve the exterior
Dirichlet scattering problem for the Helmholtz equation, we employ an
integral equation formulation whose kernel is given by the combined
field potential~$K$:
\begin{equation}
K(\bx,\by) = (\bn(\by) \cdot \nabla_{\by}G(\bx,\by)) - ik G(\bx,\by).
\end{equation}
After inverting the resulting integral equation
\begin{equation}
  \frac{1}{2} \sigma + \cK_\Gamma [\sigma] = f, \qquad \text{on } \Gamma,
\end{equation}
the solution $u$ to the boundary value problem is given as
\begin{equation}
  u(\bx) = \cK_\Gamma [\sigma](\bx), \qquad \text{for } \bx \in \bbR^3
  \setminus \Omega.
\end{equation}
Our goal in the following is to justify the method of proxy
compression, i.e. to show that the appropriate row or column spaces of
submatrices~$\mtx{A}_{\sfF \sfB}$ and~$\mtx{A}_{\sfB \sfF}$, respectively,
are spanned by proxy kernel matrices which only involve the kernel~$K$
evaluated at proxy points and the relevant set of sources or targets.

\subsubsection{Outgoing skeletonization}

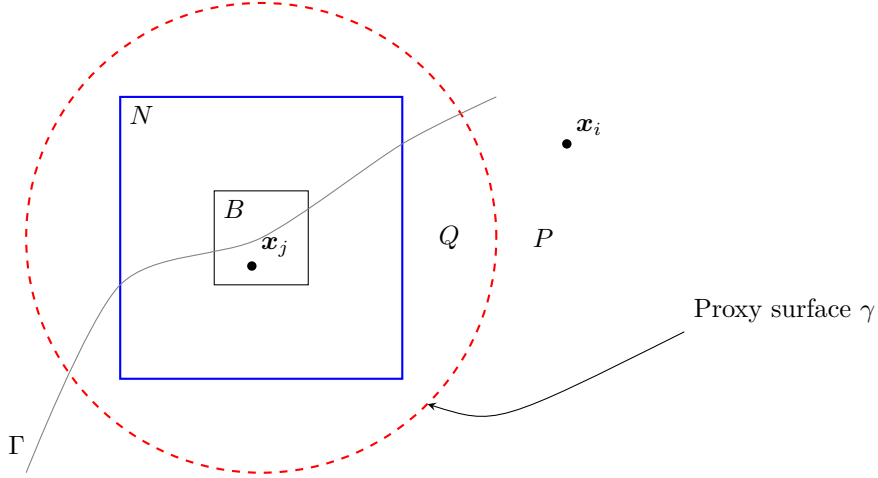
\begin{figure}[t!]
  \centering
  \begin{tikzpicture}[scale=1.25]
    \draw[blue,thick] (0,0) rectangle (3,3);
    \draw (1,1) rectangle (2,2);
    \draw[red,thick,dashed] (1.5,1.5) circle (2.5);
    \node[anchor=north west] at (1,2) {$B$};
    \fill [black] (1.4,1.2) circle (.05) node[anchor=south west] {$\bx_j$};
    \node[anchor=north west] at (0,3) {$N$};
    \node at (3.5,1.5) {$Q$};
    \node at (4.5,1.5) {$P$};
    \fill [black] (4.75,2.5) circle (.05) node[anchor=south west] {$\bx_i$};
    \draw [gray] plot [smooth] coordinates {(-1,-1) (0,1) (1.5,1.5) (3,2.5) (4,3)};
    \node[anchor=center] at (-1.1,-.7) {$\Gamma$};
    \draw[stealth-] (3.27,-.27) .. controls (4,-.5) .. (6,.5)
    node[anchor=south west] {Proxy surface $\gamma$};
  \end{tikzpicture}
  \caption{The proxy surface setup for a box~$B$. The near field is denoted
    by~$N$, and the far field has been partitioned as~$F = Q \cup P$.}
  \label{fig:proxy}
\end{figure}

To this end, for a given box~$B$, a region $D$ and its
boundary~$\gamma$ can be chosen such that $B \subset D$ and that
$\bx_{i} \in P \implies \bx_{i} \in D^c$. See Figure~\ref{fig:proxy}
for a geometric setup of the situation.  We therefore, by our earlier
assumptions, have that
$K: \mathbb{R}^{3}\setminus P \times P \to \mathbb{C}$ and that
$K(\cdot,\by)$ satisfies the underlying PDE (i.e. Helmholtz equation)
in~$\mathbb{R}^{3} \setminus P$ for each $\by \in P$.  For this
box~$B$, now consider what we will call the \emph{associated exterior
  proxy boundary value problem}:
\begin{equation}
  \label{eq:vbvp}
  \begin{aligned}
    (\Delta + k^2)  v &=  0, &\qquad &\text{in } \bbR^{3}
    \setminus D,\\
    v &= \cK_B [\tau], & &\text{on } \gamma, \\
    \lim_{r \to \infty} \, r \left(\frac{\partial v}{\partial r} - ik v
    \right) &= 0, & &
  \end{aligned}
\end{equation}
where~$\cK_B[\tau]$ is Helmholtz potential due to any function~$\tau$
supported on the piece of~$\Gamma$ contained inside box~$B$,
i.e.
\begin{equation}
  \cK_B [\tau] (\bx) = \int_{\Gamma \cap B}  K(\bx,\by) \, \tau(\by) \, da(\by),
  \qquad \text{for } \bx \in \gamma.
\end{equation}
The above boundary value problem is an exterior Dirichlet problem for
the Helmholtz equation, and can be solved using an integral equation
method by first representing~$v$ as~$v = \cK_\gamma [\mu]$ for some unknown
density~$\mu$ defined along~$\gamma$. (Note that in practice,~$\gamma$ is
usually taken to be the surface of a sphere, but there is no mathematical reason
why this must be the case.)

The solution to the boundary value problem~\eqref{eq:vbvp} is unique and can be
formally obtained, as was for~$u$ in the introduction, as
\begin{equation}
  \label{eq:vsol}
  v = \cK_\gamma \left( \frac{1}{2}
    \cI + \cK_{\gamma\gamma} \right)^{-1} \cK_B [\tau],
\end{equation}
where it is understood that the inverse operator in the middle is a
map from~$\gamma$ to~$\gamma$, i.e. that~$\cK_{\gamma \gamma}$ is a
map from~$\gamma$ to~$\gamma$ and is interpreted in the proper
principal value sense.
Finally, consider a discretization and~$L^2$ embedding~\cite{bremer_2012}
of the above form of~$v$ which maps
sources at a collection of~$\bx_j \in B$
with strengths~$t_j = \tau(\bx_j)$ to their
potentials~$v_i$ at target locations~$\bx_i \in P$:
\begin{equation}
  \begin{aligned}
    \sqrt{\mtx{D}_\sfP} \vct{v} &\approx \sqrt{\mtx{D}_\sfP} {\mtx{K}}_{\sfP \gamma} \sqrt{\mtx{D}_\gamma} \lp \mtx{I}/2 +
    \sqrt{\mtx{D}_\gamma}
    {\mtx{K}}_{\gamma \gamma} \odot \mtx{W}_{\gamma\gamma} \sqrt{\mtx{D}^{-1}_\gamma}
    \rp^{-1}
    \sqrt{\mtx{D}_\gamma}
    \mtx{K}_{\gamma \sfB} \sqrt{\mtx{D}_{\sfB}} \lp \sqrt{\mtx{D}_{\sfB}}
    \vct{t} \rp \\
    &= \sqrt{\mtx{D}_\sfP } \mtx{M}_{\sfP \gamma}
    \mtx{K}_{\gamma \sfB} \sqrt{\mtx{D}_{\sfB}} \lp \sqrt{\mtx{D}_{\sfB}}
    \vct{t}\rp,
  \end{aligned}
\end{equation}
where~$\mtx{K}_{\gamma\gamma} \odot \mtx{W}_{\gamma\gamma}$ denotes the
elementwise Hadamard product of the kernel matrix~$\mtx{K}_{\gamma\gamma}$ with
a matrix~$\mtx{W}_{\gamma\gamma}$ of quadrature corrections, the
matrices~$\mtx{D}_\ell$ contain smooth quadrature corrections along their
diagonal, and
\begin{equation}
\mtx{M}_{\sfP \gamma} = {\mtx{K}}_{\sfP \gamma} \sqrt{\mtx{D}_\gamma} \lp \mtx{I}/2 +
    \sqrt{\mtx{D}_\gamma}
    {\mtx{K}}_{\gamma \gamma} \odot \mtx{W}_{\gamma\gamma} \sqrt{\mtx{D}^{-1}_\gamma}
    \rp^{-1}
    \sqrt{\mtx{D}_\gamma}.
\end{equation}
Upon further inspection, however, assuming that the discretization of
the above integral equation along~$\gamma$ was suitably accurate and that the
source and target locations~$\bx_j$ and~$\bx_i$ were chosen to be the same as in
the original discretization of~$\Gamma$, we have that it must be true that
\begin{equation}
  \label{eq:mk}
  \mtx{K}_{\sfP \sfB} \approx \mtx{M}_{\sfP \gamma} \mtx{K}_{\gamma \sfB}
\end{equation}
due to the uniqueness of the exterior Helmholtz Dirichlet problem
(i.e. if the boundary data in~\eqref{eq:vbvp} were chosen to agree with
that induced by sources located at discretization points~$\bx_j \in
B$, then the solution~$v$ at~$\bx_i$
must agree with the potential generated via
multiplication by~$\mtx{K}_{\sfP \sfB}$). Therefore, due to the
partition of the far field~$F = Q \cup P$, it must also be true
that (after a suitable permutation of rows)
\begin{equation}
  \label{eq:afb_fact}
  \mtx{A}_{\sfF \sfB} =
  \begin{bmatrix}
    \mtx{A}_{\sfQ \sfB} \\
    \mtx{A}_{\sfP \sfB}
  \end{bmatrix}
  =
  \begin{bmatrix}
    \mtx{A}_{\sfQ \sfB} \\
      \sqrt{\mtx{D}_\sfP}
  \mtx{K}_{\sfP \sfB}
  \sqrt{\mtx{D}_\sfB}
\end{bmatrix}
\approx
  \begin{bmatrix}
    \mtx{I} & \mtx{0} \\
    \mtx{0} &  \sqrt{\mtx{D}_\sfP}
  \mtx{M}_{\sfP \gamma}
\end{bmatrix}
  \begin{bmatrix}
    \mtx{A}_{\sfQ \sfB} \\
  \mtx{K}_{\gamma \sfB}
  \sqrt{\mtx{D}_\sfB}
\end{bmatrix}.
\end{equation}
A note on dimensions of the above matrices: recall that, by
assumption,~$|\sfQ| = \cO(1)$ and therefore the bulk of the discretization nodes
in~$F$ are contained in~$P$, i.e. that~$|\sfP| = \cO(n)$. This implies that
the matrix on the very right above has dimensions which
are~$\cO(n_\gamma) \times \cO(1)$, where~$n_\gamma$ denotes the number of points
used to discretize the proxy surface. Choosing~$n_\gamma$ is discussed in the
following section. Next we detail how to use the above factorization to compute
a column skeletonization of the original submatrix~$\mtx{A}_{\sfF \sfB}$.

First, a column skeletonization of the matrix on the right
in~\eqref{eq:afb_fact} is performed which yields a decomposition
of~$\sfB$ into redundant and skeleton points,
i.e.~$\sfB = \sfR \cup \sfS$:
\begin{equation}
  \label{eq:outid}
  \begin{bmatrix}
    \mtx{A}_{\sfQ \sfB} \\
  \mtx{K}_{\gamma \sfB}
  \sqrt{\mtx{D}_\sfB}
\end{bmatrix}
\approx
  \begin{bmatrix}
    \mtx{A}_{\sfQ \sfS} \\
  \mtx{K}_{\gamma \sfS}
  \sqrt{\mtx{D}_\sfS}
\end{bmatrix}
\begin{bmatrix}
  \mtx{T}_{\sfS \sfR} & \mtx{I}
\end{bmatrix}
\mtx{P}_{\sfB}
\end{equation}
where~$\mtx{P}_{\sfB}$ is a permutation of the matrix performing a
reordering of columns according to the split~$\sfB = \sfR \cup \sfS$, and
in our continued abuse of notation, it is implied that this is
an~$\epsilon$-accurate approximation.
Inserting~\eqref{eq:outid} into~\eqref{eq:afb_fact}, we have
\begin{equation}
  \label{eq:incsr}
  \begin{aligned}
  \mtx{A}_{\sfF \sfB} &\approx
  \begin{bmatrix}
    \mtx{I} & \mtx{0} \\
    \mtx{0} &  \sqrt{\mtx{D}_\sfP}
  \mtx{M}_{\sfP \gamma}
\end{bmatrix}  
  \begin{bmatrix}
    \mtx{A}_{\sfQ \sfS} \\
  \mtx{K}_{\gamma \sfS}
  \sqrt{\mtx{D}_\sfS}
\end{bmatrix}
\begin{bmatrix}
  \mtx{T}_{\sfS \sfR} & \mtx{I}
\end{bmatrix}
\mtx{P}_{\sfB} \\
&=
\begin{bmatrix}
    \mtx{A}_{\sfQ \sfS} \\
\sqrt{\mtx{D}_\sfP}
  \mtx{M}_{\sfP \gamma} \mtx{K}_{\gamma \sfS}
  \sqrt{\mtx{D}_\sfS}
\end{bmatrix}
\begin{bmatrix}
  \mtx{T}_{\sfS \sfR} & \mtx{I}
\end{bmatrix}
\mtx{P}_{\sfB}.
\end{aligned}
\end{equation}
Again, due to the uniqueness of the exterior Helmholtz
boundary value problem in~\eqref{eq:vbvp}, it must be true
that~$\mtx{M}_{\sfP \gamma} \mtx{K}_{\gamma \sfS} \approx  \mtx{K}_{\sfP \sfS}$,
and therefore we have computed a low-rank approximation to the
original matrix~$\mtx{A}_{\sfF \sfB}$ since
\begin{equation}
  \begin{aligned}
\mtx{A}_{\sfF \sfB} &=
\begin{bmatrix}
    \mtx{A}_{\sfQ \sfS} \\
\sqrt{\mtx{D}_\sfP}
  \mtx{K}_{\sfP \sfS}
  \sqrt{\mtx{D}_\sfS}
\end{bmatrix}
\begin{bmatrix}
  \mtx{T}_{\sfS \sfR} & \mtx{I}
\end{bmatrix}
\mtx{P}_{\sfB} \\
&= \begin{bmatrix}
    \mtx{A}_{\sfQ \sfS} \\
  \mtx{A}_{\sfP \sfS}
\end{bmatrix}
\begin{bmatrix}
  \mtx{T}_{\sfS \sfR} & \mtx{I}
\end{bmatrix}
\mtx{P}_{\sfB} \\
&= \mtx{A}_{\sfF \sfS}
\begin{bmatrix}
  \mtx{T}_{\sfS \sfR} & \mtx{I}
\end{bmatrix}
\mtx{P}_{\sfB}.
\end{aligned}
\end{equation}

\begin{remark}
  In the publicly available software associated
  with~\cite{minden-str_scel-2017}, the authors use an alternative method to
  obtain interpolative decompositions of matrices such as~$\mtx{A}_{\sfF \sfB}$
  via proxy compression. It appears that in practice, their method doesn't seem
  to affect the accuracy of the computed solution on simple geometries for
  requested tolerances up to~$10^{-8}$. Our discussion above is meant to be a
  fully self-contained treatment of proxy compression, as was implemented in our
  solver, based on first-principles considerations.
\end{remark}

\subsubsection{Incoming skeletonization}

The above algorithm for compressing~$\mtx{A}_{\sfF \sfB}$ can immediately be
extended to the row compression of the dual block~$\mtx{A}_{\sfB \sfF}$ by
considering instead the \emph{associated interior proxy boundary value problem},
analogous to~\eqref{eq:vbvp}:
\begin{equation}
  \label{eq:aipbvp}
  \begin{aligned}
    (\Delta + k^2)  w &=  0, &\qquad &\text{in } D,\\
    w &= \cK_{P} [\tau], & &\text{on } \gamma,
  \end{aligned}
\end{equation}
where~$\cK_{P}[\tau]$ is Helmholtz potential due to any function~$\tau$
supported on the piece of~$\Gamma$ contained outside of~$B$, its near field~$N$, and its associated quadrature-corrected far field~$Q$,
i.e.
\begin{equation}
  \cK_{P} [\tau] (\bx) = \int_{\Gamma \cap P}  K(\bx,\by) \, \tau(\by) \, da(\by),
  \qquad \text{for } \bx \in \gamma.
\end{equation}
The above boundary value problem is an interior Dirichlet problem for the
Helmholtz equation, and can similarly be solved using an integral equation
method with one additional caveat: since it is an interior problem for the
Helmholtz equation, the solution may not be unique. That is to say, if~$-k^{2}$
is an interior Dirichlet eigenvalue for the Laplace operator in~$D$, then the
solution to the above boundary value problem is not unique. Fortunately, the
method of proxy compression does \emph{not} require uniqueness to this problem
in order to work, merely that all solutions~$w$ are spanned by some suitable
subspace of solutions to the Helmholtz equation. We now provide details to this
end {\note using an especially simple representation of the solution via a
  single layer potential, while simultaneously augmenting the right-hand side of
  the boundary value problem to be a linear combination of Dirichlet and Neumann
  data, which is readily available in most scattering problems.}

Let us begin by representing the solution~$w$ merely as a single-layer density
along~$\gamma$,
\begin{equation}
  \begin{aligned}
    w(\bx) &= \int_{\gamma} G(\bx,\by) \, \mu(\by) \, da(\by) \\
    &= \cS_{\gamma} [\mu] (\bx),
  \end{aligned}
\end{equation}
and denote the boundary data by~$f$: $f = \cK_{P} [\tau]$. Using this
representation to derive an integral equation for~$\eqref{eq:aipbvp}$ leads to a
first-kind integral equation which is not invertible if~$-k^{2}$ is an interior
eigenvalue of the Laplace operator on~$D$. However, if we also assume that we
can evaluate the normal derivative of~$f$ along~$\gamma$,
\begin{equation}
  f' = \frac{\partial f}{\partial n} = \frac{\partial}{\partial n} \cK_{P} [\tau],
\end{equation}
then we can linearly combine these boundary conditions,~$ik w + w' = ik f + f'$, obtaining the second-kind integral equation along~$\gamma$
\begin{equation}
  \lp \frac{1}{2} \cI + \cS'_{\gamma \gamma} - ik \cS_{\gamma \gamma} \rp [\mu]
  = -ik f + f'.
\end{equation}
It is well-known, and easy to verify, that the above integral operator is the
transpose (i.e. non-conjugate adjoint) of the integral operator obtained when
using the combined field representation to solve the exterior Dirichlet problem
for the Helmholtz equation, as seen in~\eqref{eq:vsol}, and therefore is
uniquely invertible~\cite{colton_kress} \emph{for all}~$k$. (The operator~$\cS_{\gamma \gamma}$ is
symmetric, and the double layer~$\cD_{\gamma \gamma}$ can be obtained
from~$\cS'_{\gamma\gamma}$ merely by switching arguments.)
 Let us denote this operator by
\begin{equation}
  \cK_{\gamma\gamma}^{T} = \frac{1}{2} \cI + \cS'_{\gamma \gamma} - ik \cS_{\gamma \gamma}.
\end{equation}
The solution to the associated interior proxy boundary value problem is then
given by
\begin{equation}
  \label{eq:wsol}
  w = \cS_\gamma \left( \frac{1}{2}\cI + \cK_{\gamma\gamma}^{T} \right)^{-1} \left[
    -ik f + f' \right].
\end{equation}
The above shows that a \emph{complete} representation of Helmholtz potentials
inside~$D$ is given in terms of a single layer potential along~$\gamma$; in
particular, this includes Helmholtz potentials which are Dirichlet
eigenfunctions of the Laplace operator. {\note Note that uniquely
  determining the solution requires the use of both~$f$ and $f'$, its normal
  derivative along~$\gamma$. This is \emph{not} usually the case for standard
  interior Dirichlet boundary value problems.} In what follows, let us set
$\tilde f = -ik f + f'$.

As in the previous section, consider a discretization and~$L^2$ embedding
of the above form of~$w$ which maps
sources at a collection of~$\bx_j \in P$
with strengths~$t_j = \tau(\bx_j)$ to their
potentials~$w_i$ at target locations~$\bx_i \in B$:
\begin{equation}
  \begin{aligned}
    \sqrt{\mtx{D}_\sfB} \vct{w} &\approx \sqrt{\mtx{D}_\sfB} {\mtx{S}}_{\sfB \gamma}
\sqrt{\mtx{D}_\gamma} \lp \mtx{I}/2 +
    \sqrt{\mtx{D}_\gamma}
    {\mtx{K}^{T}}_{\gamma \gamma} \odot \mtx{W}_{\gamma\gamma} \sqrt{\mtx{D}^{-1}_\gamma}
    \rp^{-1}
    \sqrt{\mtx{D}_\gamma}
    \tilde{ \mtx{K}}_{\gamma \sfP} \sqrt{\mtx{D}_{\sfP}} \lp \sqrt{\mtx{D}_{\sfP}}
    \vct{t} \rp \\
    &= \sqrt{\mtx{D}_\sfB } \mtx{S}_{\sfB \gamma}
    \mtx{M}_{\gamma \sfP} \sqrt{\mtx{D}_{\sfP}} \lp \sqrt{\mtx{D}_{\sfP}}
    \vct{t}\rp,
  \end{aligned}
\end{equation}
where this time
\begin{equation}
  \mtx{M}_{\gamma \sfP} =
  \sqrt{\mtx{D}_\gamma} \lp \mtx{I}/2 +
  \sqrt{\mtx{D}_\gamma}
  {\mtx{K}^{T}}_{\gamma \gamma} \odot \mtx{W}_{\gamma\gamma} \sqrt{\mtx{D}^{-1}_\gamma}
  \rp^{-1}
  \sqrt{\mtx{D}_\gamma}
  \tilde{\mtx{K}}_{\gamma \sfP},
\end{equation}
$\mtx{S}$ denotes the discretization of the single layer operator~$\cS$,
and~$\tilde{\mtx{K}}$ denotes the discretization of~\mbox{$-ik\cK + \cK'$}.
As for the
associated exterior proxy boundary value problem, assuming again that the
discretization of the above integral equation along~$\gamma$ was suitably
accurate and that the source and target locations~$\bx_j$ and~$\bx_i$ were
chosen to be the same as in the original discretization of~$\Gamma$, we have
that it must be true that there exists some matrix~$\tilde{\mtx{M}}_{\gamma \sfP}$
such that
\begin{equation}
  \label{eq:mkinc}
  \mtx{K}_{\sfB \sfP} \approx \mtx{S}_{\sfB \gamma} \tilde{\mtx{M}}_{\gamma \sfP}
\end{equation}
due to the uniqueness of the interior Helmholtz problem when both Dirichlet and Neumann data are available in the form~$-ikf + f'$.
 (In particular, it was just shown above that the range
of~$\mtx{K}_{\sfB \sfP}$ is contained in the range of~$\mtx{S}_{\sfB \gamma}$.)

Using the above factorization, the matrix~$\mtx{A}_{\sfB \sfF}$ can
be explicitly factored into the form:
\begin{equation}
  \label{eq:abf_fact}
  \begin{aligned}
  \mtx{A}_{\sfB \sfF} &=
  \begin{bmatrix}
    \mtx{A}_{\sfB \sfQ} &
    \mtx{A}_{\sfB \sfP}
  \end{bmatrix} \\
  &=
  \begin{bmatrix}
    \mtx{A}_{\sfB \sfQ} &
      \sqrt{\mtx{D}_\sfB}
  \mtx{K}_{\sfB \sfP}
  \sqrt{\mtx{D}_\sfP}
\end{bmatrix} \\
&\approx
  \begin{bmatrix}
    \mtx{A}_{\sfB \sfQ} &
  \sqrt{\mtx{D}_\sfB}
  \mtx{S}_{\sfB \gamma}
\end{bmatrix}
  \begin{bmatrix}
    \mtx{I} & \mtx{0} \\
    \mtx{0} &  \tilde{\mtx{M}}_{\gamma \sfP} \sqrt{\mtx{D}_\sfP}
  \end{bmatrix}.
\end{aligned}
\end{equation}
Next, a row skeletonization of the matrix
\begin{equation}
  \begin{bmatrix}
    \mtx{A}_{\sfB \sfQ} &
  \sqrt{\mtx{D}_\sfB}
  \mtx{S}_{\sfB \gamma}
\end{bmatrix}
\end{equation}
can be computed, splitting the row indices into redundant and skeleton
points,~$\sfB = \sfR \cup \sfS$:
\begin{equation}
  \begin{bmatrix}
    \mtx{A}_{\sfB \sfQ} &
  \sqrt{\mtx{D}_\sfB}
  \mtx{S}_{\sfB \gamma}
  \end{bmatrix}
  \approx \mtx{P}_{\sfB}
  \begin{bmatrix}
    \mtx{T}_{\sfR \sfS} \\
    \mtx{I}
  \end{bmatrix}
  \begin{bmatrix}
    \mtx{A}_{\sfS \sfQ} &
  \sqrt{\mtx{D}_\sfS}
  \mtx{S}_{\sfS \gamma}
  \end{bmatrix},
\end{equation}
where~$\mtx{P}_{\sfB}$ is a suitable permutation matrix which appropriately
reorders the skeleton and redundant points. Finally, it can be shown that the
matrix~$\mtx{A}_{\sfB \sfF}$ can be written (i.e. compressed) as:
\begin{equation}
  \label{eq:abf_fact2}
  \begin{aligned}
    \mtx{A}_{\sfB \sfF} &\approx
\mtx{P}_{\sfB}
  \begin{bmatrix}
    \mtx{T}_{\sfR \sfS} \\
    \mtx{I}
  \end{bmatrix}
  \begin{bmatrix}
    \mtx{A}_{\sfS \sfQ} &
  \sqrt{\mtx{D}_\sfS}
  \mtx{S}_{\sfS \gamma}
  \end{bmatrix}
  \begin{bmatrix}
    \mtx{I} & \mtx{0} \\
    \mtx{0} &  \tilde{\mtx{M}}_{\gamma \sfP} \sqrt{\mtx{D}_\sfP}
  \end{bmatrix} \\
&\approx \mtx{P}_{\sfB}
  \begin{bmatrix}
    \mtx{T}_{\sfR \sfS} \\
    \mtx{I}
  \end{bmatrix}
  \begin{bmatrix}
    \mtx{A}_{\sfS \sfQ} & \mtx{A}_{\sfS \sfP}
  \end{bmatrix} \\
&= \mtx{P}_{\sfB}
  \begin{bmatrix}
    \mtx{T}_{\sfR \sfS} \\
    \mtx{I}
  \end{bmatrix}
  \begin{bmatrix}
    \mtx{A}_{\sfS \sfF}
  \end{bmatrix},
  \end{aligned}
\end{equation}
where we invoked the fact that, up to discretization
error,~$\mtx{S}_{\sfS \gamma} \tilde{\mtx{M}}_{\gamma \sfP} \approx \mtx{K}_{\sfS \sfP}$.
We refer to this procedure as calculating the
\emph{incoming} skeleton nodes for a box~$B$.

{\note
Lastly,
we note again that the matrix~$\mtx{M}_{\gamma \sfP}$ in~\eqref{eq:abf_fact2}
  is only ever used implicitly --  this matrix is
never explicitly formed, it is only the matrices~$\mtx{A}_{\sfB \sfQ}$
and~$\sqrt{\mtx{D}_\sfB} \mtx{S}_{\sfB \gamma}$ which need to be formed and
compressed. Put another way, $w = \cS_\gamma[\mu]$ is a complete representation
of solutions to the interior boundary value problem~\eqref{eq:aipbvp} despite
the boundary value problem possibly not having a unique solution~\cite{colton_kress}.
}

\subsubsection{Simultaneous skeletonization}

Lastly, as mentioned in Section~\ref{subsec:proxy}, it is sometimes
useful to choose identical incoming and outgoing skeleton points for a
particular box. If this is desired, then a single interpolative
decomposition can be performed:
\begin{equation}
  \begin{bmatrix}
    \mtx{A}_{\sfQ \sfB} \\
    \mtx{A}_{\sfB \sfQ}^T \\
  \mtx{K}_{\gamma \sfB}
  \sqrt{\mtx{D}_\sfB} \\
  \mtx{S}_{\sfB \gamma}^T
  \sqrt{\mtx{D}^T_\sfB}
\end{bmatrix}
\approx
  \begin{bmatrix}
    \mtx{A}_{\sfQ \sfS} \\
    \mtx{A}_{\sfS \sfQ}^T \\
  \mtx{K}_{\gamma \sfS}
  \sqrt{\mtx{D}_\sfS} \\
  \mtx{S}_{\sfS \gamma}^T
  \sqrt{\mtx{D}^T_\sfS}
\end{bmatrix}
\begin{bmatrix}
  \mtx{T}_{\sfS \sfR} & \mtx{I}
\end{bmatrix}
\mtx{P}_{\sfB}.
\end{equation}
It is easily verified that the matrix~$\mtx{T}_{\sfS \sfR}$ above
compresses both~$\mtx{A}_{\sfF \sfB}$ and~$\mtx{A}_{\sfB \sfF}^T$. In
practice, the rank used in the above low-rank approximation is
slightly larger than if the blocks had been compressed separately
(but, the subsequent bookkeeping and factorization is slightly simpler
and uses less storage).

\subsection{Proxy surface discretization}

When using the proxy surface compression technique, the
matrix~$\mtx{K}_{\gamma \sfB}$ in~\eqref{eq:mk} must accurately capture
the subspace of potentials supported on the proxy surface, and
likewise,~$\mtx{S}_{\sfB \gamma}$ must accurately capture the subspace
of potentials \emph{induced} by densities on, or inside, the proxy
surface~$\gamma$. We therefore must take some care when choosing how
to select points~$\by_\ell$ on~$\gamma$.

{\note When the underlying Green's function for the problem is non-oscillatory
  (e.g. $G(r) = 1/r$), or when the box~$B$ that the proxy surface encloses is
  less than a wavelength in size when the kernel is oscillatory, the rank
  of~$\mtx{A}_{\sfF \sfB}$ is bounded by a constant independent of the size of
  the box~$B$.} For oscillatory problems, it turns out that the rank of
$\mtx{A}_{\sfF \sfB}$ tends to grow with the size of the box~$B$ measured in
terms of wavelengths. This growth in rank can be attributed to the oscillatory
nature of the Green's function $G$, which, in turn, necessitates an increased
number of points to sample the proxy surface for accurately computing the
implied integrals and subsequent compressed representation. One can use standard
multipole estimates for appropriately choosing~$n_{\gamma}$ to enable accurate
compression of the far field blocks~\cite{wideband3d,greengard-1997}.

{\note As the FMM-LU algorithm proceeds up the octree data structure, the
  region~$B$ which contains the active degrees of freedom grows, as does the
  size (and number of discretization points) in the near field of~$B$. } In
particular, in the oscillatory regime, the rank of the interaction of a box~$B$
with its far field~$F$ tends to scale quadratically with the diameter of the
box~$B$ measured in wavelengths. {\note When the size of the problem as measured
  in wavelengths is fixed, in the refinement limit $n\to \infty$, the number of
  points required on the proxy surface for a given accuracy~$n_{\gamma}$ is
  still a constant, and typically $\ll n$. In this fixed-wavenumber regime,
  refinement does not affect the overall complexity of the method.}
However, in the regime where the surface~$\Gamma$, and therefore the proxy
surfaces~$\gamma$, are sampled at fixed number of points per wavelength, then
the cost of the factorization typically grows like~$\cO(n^2)$
since~$k^{2} \propto n$. This is the expected behavior of fast direct solvers
which are based on non-directional low-rank approximation, as demonstrated
in~\cite{bremer-2015, gopal2020ls, martinsson-book}. The fact is that these
interactions cease to be low-rank in the highly oscillatory regime.

\begin{remark}
  When finding the skeleton and redundant points in a box, the far interaction
  matrices and the blocks affected by Schur complements are simultaneously
  compressed. Due to the presence of near blocks which contain the Schur
  complements, in practice, for small wave numbers, it is possible to use fewer
  points than predicted by standard multipole estimates.
\end{remark}

\subsection{Far field partitioning}
\label{sec:farsplit}

Recall that a necessary criterion for determining the split $F = Q
\cup P$ for a box $B$ is that the interactions corresponding to
points $\bx_{j} \in B$ and the points $\bx_{i} \in P$ are of the
form $\sqrt{w_{i}} K(\bx_{i}, \bx_{j}) \sqrt{w_{j}}$. There are two
mechanisms by which the condition would be violated:
\begin{enumerate}
  \item If $\bx_{i}\in \text{Near}(\bx_{j})$ from the quadrature perspective,
        then the original matrix entry has target dependent weights and is of
        the form $\sqrt{w_i} K(\bx_{i}, \bx_{j}) w_{ij}/\sqrt{w_j}$, or
  \item If the matrix entries have been updated due to Schur complements arising
        from a previous application of the strong skeletonization procedure to a
        different box.
\end{enumerate}

Depending on the subdivision criterion used for generating the octree,
the near interactions associated with the locally-corrected quadrature
may spill over into the far field region of the box, i.e. for a point
$\bx_{j} \in B$, there may exist $\bx_{i} \in \text{Near}(\bx_{j})
\cap F$, where $F$ is the far field region associated with box
$B$. In this situation, the partition $F = Q \cup P$ can be
modified to include all such near quadrature corrections of point that
spill over into the far region of the corresponding box, i.e. if
$\bx_{i} \in \text{Near}(\bx_{j}) \cap F$ for any $\bx_{j} \in B$
then $\bx_{i}\in Q$ even if it would not have been affected by
previously computed Schur complement blocks. In practice, this
situation is rarely encountered and typically happens for at most
$\cO(1)$ boxes even in multiscale geometries.

\section{Numerical experiments}
\label{sec:numerical}

\begin{figure}[t]
  \centering
    \includegraphics[width=.35\linewidth]{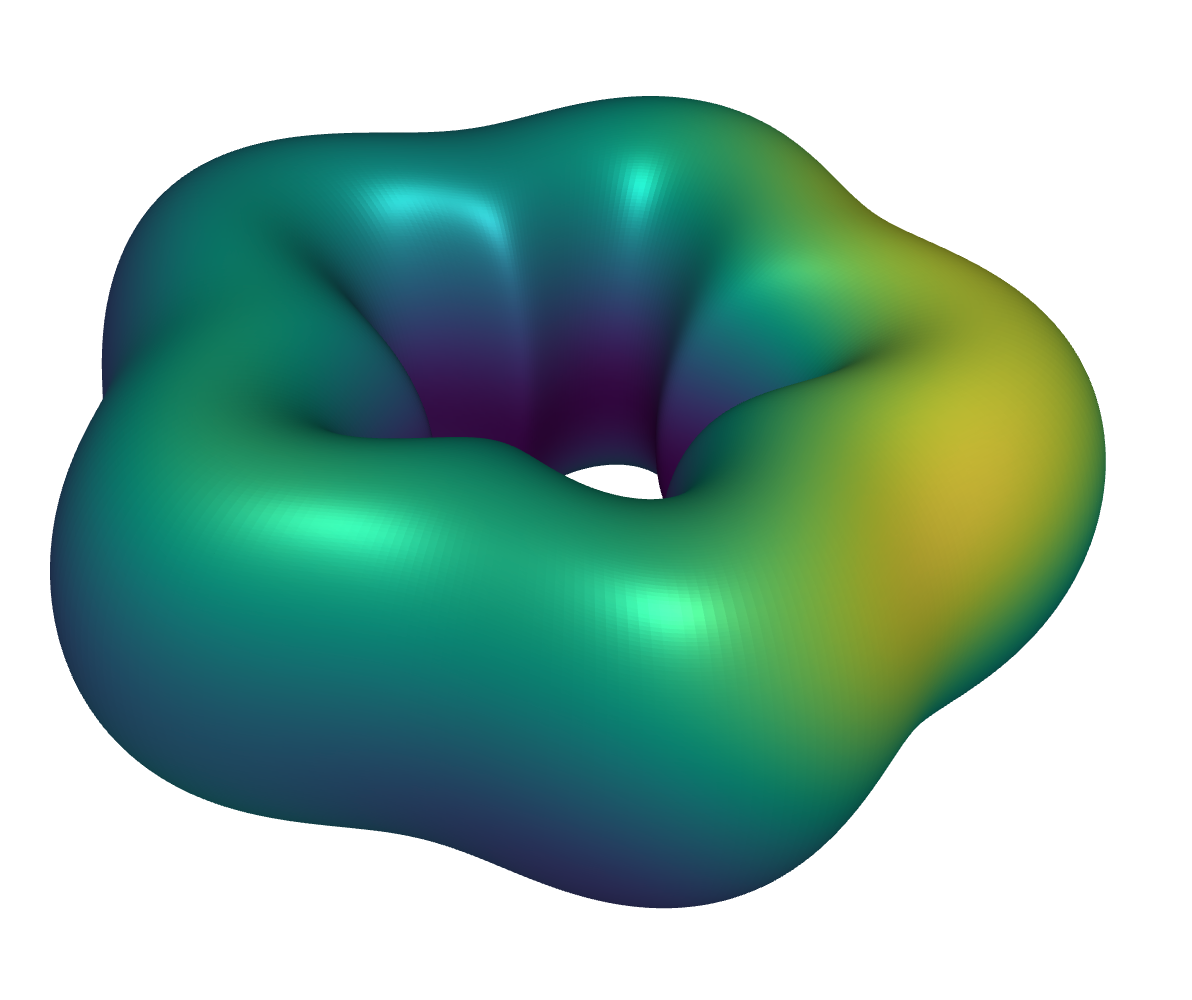}
  \caption{A wiggly torus geometry with parametrization $\bX$ defined in Equation~\cref{eq:wtorus}}
  \label{fig:wtorus}
\end{figure}

In this section, we illustrate the performance of our approach.  For
examples in~\cref{subsec:num-conv,subsec:num-freq}, we consider a
wiggly torus as the geometry where the boundary $\Gamma$ is
parameterized by~\mbox{$\bX :[0,2\pi]^2 \to \Gamma$} as
\begin{equation}
\label{eq:wtorus}
\bX(u,v) = \begin{bmatrix}
1.2 \cdot \left(2 + \cos{(v)} + 0.25 \cos{(5u)} \right) \cos{(u)} \\
\left(2 + \cos{(v)} + 0.25 \cos{(5u)} \right) \sin{(u)} \\
1.7 \cdot \sin{(v})
\end{bmatrix} \, ,
\end{equation}
see~\cref{fig:wtorus}.  In~\cref{sec:plane}, the geometry is a
multiscale plane where the ratio of the largest triangle to the
smallest one is $\cO(10^{3})$. For all three examples, we consider the
solution of the exterior Dirichlet problem~\cref{eq:extdir} using the
combined field representation
\begin{equation}
u(\bx) = \cD_{k} [\sigma](\bx) - ik \cS[\sigma](\bx) \,, \quad \bx
\in \mathbb{R}^{3} \setminus \Omega.
\end{equation}
Imposing Dirichlet boundary conditions on this expression results in the
following integral equation for the unknown density $\sigma$:
\begin{equation}
\label{eq:num-inteq}
\frac{\sigma(\bx)}{2} + \cD_{k}[\sigma](\bx) - ik \cS[\sigma](\bx)
= f(\bx) \, , \quad \bx \in \Gamma\, .
\end{equation}
Here $f$ is the prescribed Dirichlet data,~$k$ is the Helmholtz
wavenumber, and~$\cD$,~$\cS$ are interpreted using their
on-surface principal value senses.  Let $\lambda = 2 \pi/k$ denote the
corresponding wavelength.

Suppose that~\cref{eq:num-inteq} is discretized using the Nystr\"{o}m method
where the layer potentials are evaluated using the locally corrected quadrature
methods discussed in~\cite{greengard2021fast}. We assume that the surface
$\Gamma$ is given by the disjoint union of patches $\Gamma_{j}$,
$\Gamma = \cup_{j=1}^{\Npatches} \Gamma_{j}$, where each patch is parametrized
by a non-degenerate chart $\bX^{j}: T_{0} \to \Gamma_{j}$, with
$T_{0} = \{ (u,v) : u>0, v>0, u+v\leq 1 \}$ being the standard right simplex. We
further assume that the charts $\bX^{j}$, their derivative information, the
density $\sigma|_{\Gamma_{j}}$ and the data $f|_{\Gamma_{j}}$ are discretized
using order $p$ Vioreanu-Rokhlin nodes on
$T_{0}$~\cite{vioreanu_2014,greengard2021fast}.

Let $n = \Npatches p (p+1)/2$ be the total number of discretization points. Let
$\varepsilon$ denote the tolerance used for computing the factorization (i.e.
the tolerance requested in each interpolative decomposition performed). For a
patch $\Gamma_{j}$, let $R_{j}$ denote the diameter of the patch and let
$\ppw_{j} = R_{j} k/\pi$ denote the effective sampling rate on $\Gamma_{j}$
measured in points per wavelength. Let $\ppwmin = \min_{j} \ppw_{j}$. Let
$t_{f},t_{s}$ denote the factorization time and the solve time for computing
FMM-LU factorization including quadrature corrections, respectively. Let
$\tinit$ denote the problem setup time which is dominated by the time taken to
generate the quadrature corrections. Let $m_{f}$ denote the total memory required
 for storing the factorization in~\texttt{GB}, and
$S_{f}=n/t_{f}, S_{s}=n/t_{s}$, and $\sinit=n/\tinit$ denote the corresponding
speeds for the different tasks measured in points processed per second. Lastly,
let $n_{0}$ denote the size of the system matrix inverted directly at the root
level.

The fast direct solver was implemented in Matlab as a fork
of~\url{github.com/victorminden/strong-skel}, and relies on interfaces contained
in the~\texttt{fmm3dbie} library~\cite{greengard2021fast}, freely available
at~\url{github.com/fastalgorithms/fmm3dbie}, and the~\texttt{FLAM} library,
freely available at~\url{github.com/klho/FLAM}. The code used for the following
experiments and timings can be obtained
at~\url{github.com/fastalgorithms/strong-skel}; a python implementation
can be obtained
at~\url{gitlab.com/fastalgorithms/adaptive-direct-h2}. The experiments were run
on a single core of a linux workstation with $192$~GB RAM. The current
implementation is to demonstrate the accuracy of our approach and the complexity
scaling with respect to the total number of discretization points for the
various parts of the algorithm. An improved solver with better computational
performance will be released soon.

\subsection{Accuracy and convergence}
\label{subsec:num-conv}

To test the accuracy of the solver, suppose that the Dirichlet data is the
potential due to a collection of $50$ random point sources located in the
interior of~$\Omega$ given by
\begin{equation}
f(\bx) = \sum_{j=1}^{50} q_{j} \frac{e^{ik|\bx -
    \bx_{j}|}}{|\bx-\bx_{j}|} \, .
\end{equation}
In this setting, due to uniqueness, the exact solution in the exterior
is given by the right hand side of the above expression evaluated at
$\bx \in \Omega^c$. Let $\ucomp$ denote the numerically computed
solution and let $\varepsilon_{a}$ denote the relative $L^2$ error at
$50$ random target locations in the exterior denoted by $\bt_{j}$,
$j=1,2,\ldots 50$, given by
\begin{equation}
\varepsilon_{a} = \frac{\sqrt{\sum_{j=1}^{50} |\ucomp(\bt_{j}) -
    u(\bt_{j})|^2}}{\sqrt{\int_{\Gamma} |\sigma(\bx)|^2 \, {\textrm
      d}a (\bx) }} \, .
\end{equation}

In~\cref{tab:conv}, we tabulate $t_{f}$, $t_{s}$, $t_{q}$ and $m_{f}$ for $p=6$
and $\varepsilon= 5 \times 10^{-7}$. The data are plotted in Figure~\ref{fig:graphs}.
For these experiments $k=0.97$ is chosen
such that the torus is contained in a bounding box of size
$1.2 \lambda \times 1 \lambda \times 0.5 \lambda$, and
$\varepsilon = 5\times10^{-7}$. We observe the expected linear scaling in
$t_{q}$. However, for $t_{f},t_{s}$, and $m_{f}$ we observe sub-linear scaling
due to the wavenumber being held fixed while~$n$ is increased. The additional
degrees of freedom introduced beyond a certain sampling in points per wavelength
are more easily compressed. We observe the expected convergence rate of
$\max{(h^{p-1}, \varepsilon)}$ for $\varepsilon_{a}$, see Figure~\ref{fig:rates}.
This is consistent with
the analysis in~\cite{Atkinson95}. However, for $p=4$, and $\Npatches=12800$, we
also do not observe the expected order of convergence to the desired tolerance.
This can be explained by the fact that the matrix entries corresponding to the
far-interactions were generated using the underlying smooth quadrature as
opposed to an appropriately oversampled smooth quadrature rule.

\afterpage{

\begin{table}[t!]
\centering
\begin{tabular}{ccccccccc}\toprule
$p$ & $\Npatches$ & $n$ & $k$ & $t_{f} (s)$ & $t_{s}$ (s) & $t_{q}$ (s) & $m_{f}$ (GB) & $\varepsilon_{a}$  \\ \midrule
4 & 98 & 980 & 0.97 & 0.6 & 0.03 & 1.2 & 0.02 & $3.7 \times 10^{-3}$ \\
4 & 200 & 2000 & 0.97 & 3.2 & 0.04 & 2.2 & 0.07 & $5.3 \times 10^{-4}$ \\
4 & 800 & 8000 & 0.97 &52.5 & 0.1 & 7.6 & 0.6 & $8.6 \times 10^{-5}$ \\
4 & 3200 & 32000 & 0.97 &273.4 & 0.3 & 28.5 & 2.3 & $9.8 \times 10^{-6}$ \\
4 & 7200 & 72000 & 0.97 & 540.6 & 0.8 & 72.4 & 4.3 & $2.9 \times 10^{-6}$ \\
4 & 12800 & 128000 & 0.97 & 861.2 & 1.2 & 123.2 & 7.1 & $1.2 \times 10^{-6}$ \\ \midrule
6 & 98 & 2058 & 0.97 &3.6 & 0.03 & 4.3 & 0.07& $3.8 \times 10^{-4}$\\
6 & 200 & 4200 & 0.97 &17.6 & 0.06 & 8.0 & 0.2 & $1.7 \times 10^{-5}$ \\
6 & 800 & 16800 & 0.97 &159.3 & 0.2 & 28.1 & 1.4 & $8.9 \times 10^{-7}$ \\
6 & 3200 & 67200 & 0.97 & 624.3 & 0.8 & 113.0 & 4.5 & $2.8 \times 10^{-8}$\\
6 & 7200 & 151200 & 0.97 & 1128.6 & 1.5 & 229.9 & 8.8& $1.1 \times 10^{-8}$ \\
6 & 12800 & 268800 & 0.97 & 1667.5 & 2.4 & 402.4 & 14.3& $1.1 \times 10^{-8}$ \\ \midrule
8 & 98 & 3528 & 0.97 & 11.9 & 0.04 & 10.4 & 0.2 & $5.0 \times 10^{-5}$\\
8 & 200 & 7200 & 0.97 & 48.2 & 0.09 & 18.1 & 0.5 & $2.6 \times 10^{-6}$\\
8 & 800 & 28800 & 0.97 & 263.3 & 0.3 & 66.9 & 2.2 & $1.1 \times 10^{-8}$\\
8 & 3200 &115200 & 0.97 & 846.8 & 1.1 & 287.2 & 7.0 & $1.2 \times 10^{-8}$\\
8 & 7200 & 259200 & 0.97 & 1612.7 & 2.3 & 576.5 & 13.7 & $1.2 \times 10^{-8}$\\
8 & 12800 & 460800 & 0.97 & 2615.6 & 4.2 & 1065.4 & 23.0 & $1.1 \times 10^{-8}$\\ \bottomrule
\end{tabular}
\caption{Scaling results in $t_{f}$, $t_{s}$, $t_{q}$ and $m_{f}$ as
  $p$ and $\Npatches$ is varied. The results are for fixed $k=0.97$
  and $\varepsilon=5 \times 10^{-7}$.}
\label{tab:conv}
\end{table}

\clearpage
}

\afterpage{
\begin{figure}[t!]
  \centering
  \begin{subfigure}[t]{.45\linewidth}
    \centering
    \begin{tikzpicture}[scale=0.8]
      \begin{loglogaxis} [
        legend pos = south east,
        xlabel = {$n$},
        ylabel = {$t_{f}$ (sec)}
        ]
        \addplot[
        color = black,
        mark = square,
        style = thick,
        ]
        coordinates {
          (980,0.6) (2000,3.2) (8000,52.5) (32000,273.4) (72000,540.6) (128000,861)
        };
        \addlegendentry{$p = 4$};

        \addplot[
        color = blue,
        mark = triangle,
        style = thick,
        ]
        coordinates {
          (2058,3.6) (4200,17.6) (16800,159.3) (67200,624.3) (151200,1128.6) (268800,1667.6)
        };
        \addlegendentry{$p = 6$};

        \addplot[
        color = red,
        mark = diamond,
        style = thick,
        ]
        coordinates {
          (3528,11.9) (7200,48.2) (28800,263.3) (115200,846.8) (259200,1612.7) (460800,2615.6)
        };
        \addlegendentry{$p = 8$};

        \addplot[
        color = gray,
        style = {thick, dashed},
        ]
        coordinates {
          (1000,5) (1000000,5000)
        };
        \addlegendentry{$\mathcal O(n)$};

      \end{loglogaxis}
    \end{tikzpicture}
    \caption{Factorization time.}
  \end{subfigure} \quad
  \begin{subfigure}[t]{.45\linewidth}
    \centering
    \begin{tikzpicture}[scale=0.8]
      \begin{loglogaxis} [
        legend pos = south east,
        xlabel = {$n$},
        ylabel = {$t_{s}$ (sec)}
        ]
        \addplot[
        color = black,
        mark = square,
        style = thick,
        ]
        coordinates {
          (980,0.03) (2000,0.04) (8000,0.1) (32000,0.3) (72000,0.8) (128000,1.2)
        };
        \addlegendentry{$p = 4$};

        \addplot[
        color = blue,
        mark = triangle,
        style = thick,
        ]
        coordinates {
          (2058,0.03) (4200,0.06) (16800,0.2) (67200,0.8) (151200,1.5) (268800,2.4)
        };
        \addlegendentry{$p = 6$};

        \addplot[
        color = red,
        mark = diamond,
        style = thick,
        ]
        coordinates {
          (3528,0.04) (7200,0.09) (28800,0.3) (115200,1.1) (259200,2.3) (460800,4.2)
        };
        \addlegendentry{$p = 8$};

        \addplot[
        color = gray,
        style = {thick, dashed},
        ]
        coordinates {
          (1000,0.01) (1000000,10)
        };
        \addlegendentry{$\mathcal O(n)$};

      \end{loglogaxis}
    \end{tikzpicture}
    \caption{Solve time.}
  \end{subfigure}\\
  \vspace{\baselineskip}
  \begin{subfigure}[t]{.45\linewidth}
    \centering
    \begin{tikzpicture}[scale=0.8]
      \begin{loglogaxis} [
        legend pos = south east,
        xlabel = {$n$},
        ylabel = {$t_{q}$ (sec)}
        ]
        \addplot[
        color = black,
        mark = square,
        style = thick,
        ]
        coordinates {
          (980,1.2) (2000,2.2) (8000,7.6) (32000,28.5) (72000,72.4) (128000,123.2)
        };
        \addlegendentry{$p = 4$};

        \addplot[
        color = blue,
        mark = triangle,
        style = thick,
        ]
        coordinates {
          (2058,4.3) (4200,8.0) (16800,28.1) (67200,113) (151200,229.9) (268800,402.4)
        };
        \addlegendentry{$p = 6$};

        \addplot[
        color = red,
        mark = diamond,
        style = thick,
        ]
        coordinates {
          (3528,10.4) (7200,18.1) (28800,66.9) (115200,287.2) (259200,576.5) (460800,1065.4)
        };
        \addlegendentry{$p = 8$};

        \addplot[
        color = gray,
        style = {thick, dashed},
        ]
        coordinates {
          (1000,2) (1000000,2000)
        };
        \addlegendentry{$\mathcal O(n)$};

      \end{loglogaxis}
    \end{tikzpicture}
    \caption{Initialization time (quadrature generation).}
  \end{subfigure} \quad
  \begin{subfigure}[t]{.45\linewidth}
    \centering
    \begin{tikzpicture}[scale=0.8]
      \begin{loglogaxis} [
        legend pos = south east,
        xlabel = {$n$},
        ylabel = {$m_{f}$ (GB)}
        ]
        \addplot[
        color = black,
        mark = square,
        style = thick,
        ]
        coordinates {
          (980,0.02) (2000,0.07) (8000,0.6) (32000,2.3) (72000,4.3) (128000,7.1)
        };
        \addlegendentry{$p = 4$};

        \addplot[
        color = blue,
        mark = triangle,
        style = thick,
        ]
        coordinates {
          (2058,0.07) (4200,0.2) (16800,1.4) (67200,4.5) (151200,8.8) (268800,14.3)
        };
        \addlegendentry{$p = 6$};

        \addplot[
        color = red,
        mark = diamond,
        style = thick,
        ]
        coordinates {
          (3528,0.2) (7200,0.5) (28800,2.2) (115200,7.0) (259200,13.7) (460800,23)
        };
        \addlegendentry{$p = 8$};

        \addplot[
        color = gray,
        style = {thick, dashed},
        ]
        coordinates {
          (1000,0.05) (1000000,50)
        };
        \addlegendentry{$\mathcal O(n)$};

      \end{loglogaxis}
    \end{tikzpicture}
    \caption{Memory requirements.}
  \end{subfigure}
  \caption{Timings and memory usage for the FMM-LU solver, based on data in Table~\ref{tab:conv}.}
  \label{fig:graphs}
\end{figure}
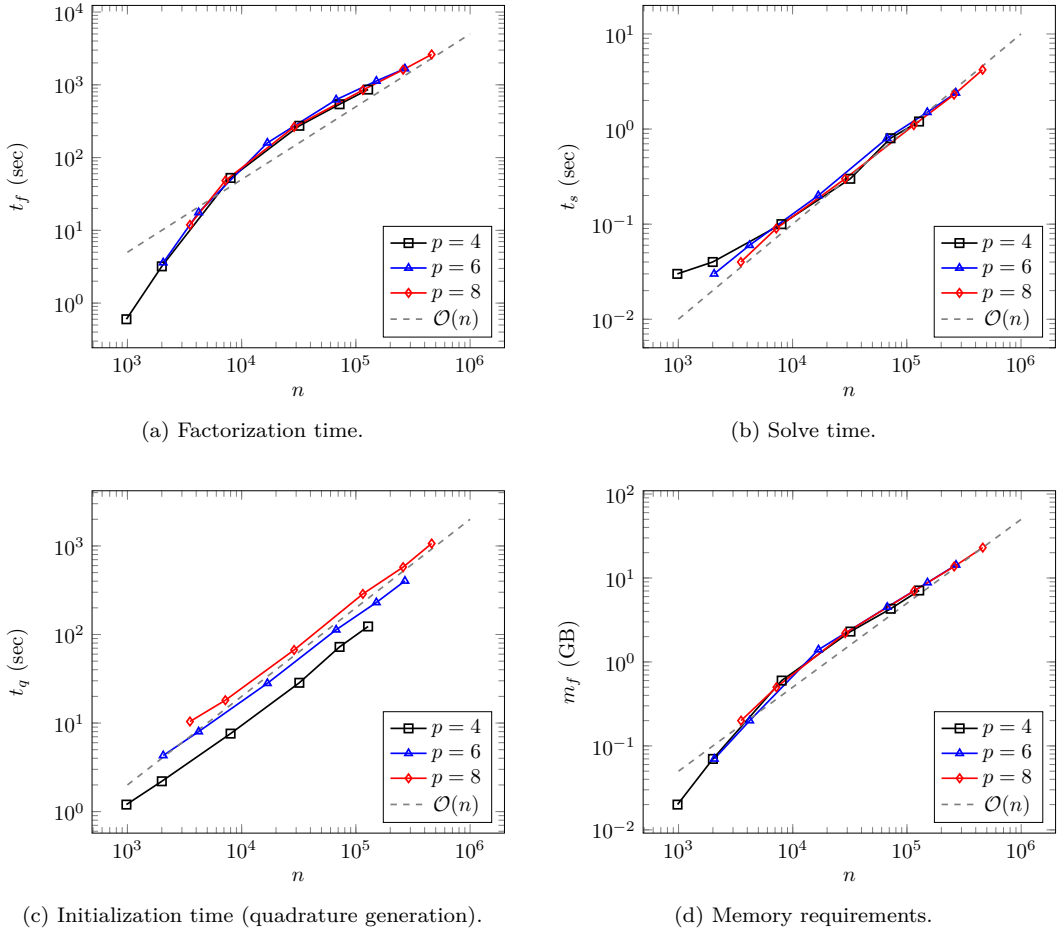

\begin{figure}[b!]
  \centering
    \begin{tikzpicture}[scale=0.8]
      \begin{loglogaxis} [
        legend pos = north east,
        xlabel = {$\Npatches$},
        ylabel = {$\varepsilon_{a}$}
        ]
        \addplot[
        color = red,
        mark = square,
        style = thick,
        ]
        coordinates {
          (98,0.0037) (200,0.00053) (800,0.000086) (3200,0.0000098) (7200,0.0000029) (12800,0.0000012)
        };
        \addlegendentry{$p = 4$};

        \addplot[
        color = green,
        mark = triangle,
        style = thick,
        ]
        coordinates {
          (98,0.00038) (200,0.000017) (800,0.00000089) (3200,0.000000028)
          (7200,0.000000011) (12800,0.000000011)
        };
        \addlegendentry{$p = 6$};

        \addplot[
        color = blue,
        mark = diamond,
        style = thick,
        ]
        coordinates {
          (98,0.00005) (200,0.0000026) (800,0.000000011) (3200,0.000000012)
          (7200,0.000000012) (12800,0.000000011)
        };
        \addlegendentry{$p = 8$};

         \addplot[
         color = red,
         style = {dashed},
         ]
         coordinates {
           (50,0.005) (12800,0.0000012)
         };
         \addlegendentry{$\mathcal O(h^{3})$};

         \addplot[
         color = green,
         style = {dashed},
         ]
         coordinates {
           (50,0.001) (3200,0.000000031)
         };
         \addlegendentry{$\mathcal O(h^{5})$};

         \addplot[
         color = blue,
         style = {dashed},
         ]
         coordinates {
           (50,0.0005) (800,0.0000000305)
         };
         \addlegendentry{$\mathcal O(h^{7})$};

      \end{loglogaxis}
    \end{tikzpicture}
    \caption{Convergence rates, based on data in Table~\ref{tab:conv}.}
  \label{fig:rates}
\end{figure}
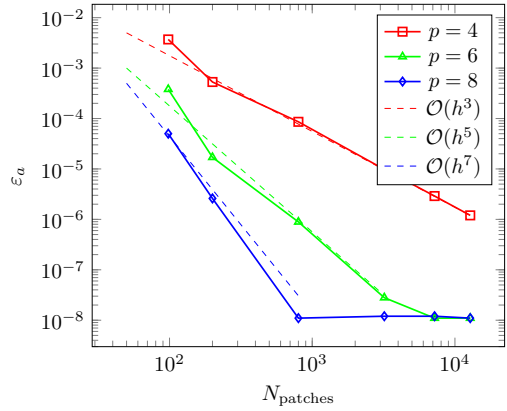

\clearpage
}

\afterpage{

\begin{table}[t!]
  \centering
  \begin{tabular}{cccccccccc}\toprule
  $p$ & $\Npatches$ & $n$ & $k$ & $t_{f} (s)$ & $t_{s}$ (s) & $t_{q}$ (s) & $m_{f}$ (GB) & $n_{0}$ & $\varepsilon_{a}$ \\ \midrule
  6 & 3200 & 67200 & 3.88 & 724.0 & 0.7 & 130.2 & 5.3 & 5700 & $1.4 \times 10^{-7}$ \\
  6 & 7200 & 151200 & 5.82 & 1511.8 & 1.9 & 260.8 & 10.8 & 7391 & $5.3 \times 10^{-8}$ \\
  6 & 12800 & 268800 & 7.76 & 2531.2 & 3.0 & 470.1 & 18.2 & 9109 & $7.3 \times 10^{-8}$\\
  6 & 28800 & 604800 & 11.64 & 5505.9 & 7.8 & 1077.3 & 38.1 & 12805 & $1.3 \times 10^{-7}$ \\
  6 & 51200 & 1075200 & 15.52 & 9756.3 & 12.1 & 1936.4 & 65.5 & 16921 & $1.5 \times 10^{-7}$ \\
  \bottomrule
  \end{tabular}
  \caption{Scaling results in $t_{f}$, $t_{s}$, $t_{q}$ and $m_{f}$ as
    $p$ and $\Npatches$ is varied. The results are for fixed
    $\ppwmin\approx 10$ and $\varepsilon=5 \times 10^{-7}$.}
  \label{tab:freq}
  \end{table}
  
  \vspace{\baselineskip}

\begin{figure}[b!]
  \centering
  \begin{subfigure}[t]{.45\linewidth}
    \centering
    \begin{tikzpicture}[scale=0.8]
      \begin{loglogaxis} [
        legend pos = south east,
        xlabel = {$n$},
        ylabel = {$t_{f}$ (sec)}
        ]
        \addplot[
        color = black,
        mark = square,
        style = thick,
        ]
        coordinates {
          (67200,724.0) (151200, 1511.8) (268800,2531.2) (604800,5505.9) 
          (1075200,9756.3)
        };
        \addlegendentry{$n \propto k^2$};

        \addplot[
        color = gray,
        style = {thick, dashed},
        ]
        coordinates {
          (67200,640) (1075200,10240)
        };
        \addlegendentry{$\mathcal O(n)$};

        \addplot[
        color = gray,
        style = {thick, dotted},
        ]
        coordinates {
          (67200,640) (1075200,25803)
        };
        \addlegendentry{$\mathcal O(n^{4/3})$};

      \end{loglogaxis}
    \end{tikzpicture}
    \caption{Factorization time.}
  \end{subfigure} \quad
  \begin{subfigure}[t]{.45\linewidth}
    \centering
    \begin{tikzpicture}[scale=0.8]
      \begin{loglogaxis} [
        legend pos = south east,
        xlabel = {$n$},
        ylabel = {$t_{s}$ (sec)}
        ]
        \addplot[
        color = black,
        mark = square,
        style = thick,
        ]
        coordinates {
          (67200,0.7) (151200, 1.9) (268800,3.0) (604800,7.8) 
          (1075200,12.1)
        };
        \addlegendentry{$n \propto k^2$};

        \addplot[
        color = gray,
        style = {thick, dashed},
        ]
        coordinates {
          (67200,.75) (1075200,12)
        };
        \addlegendentry{$\mathcal O(n)$};

        \addplot[
        color = gray,
        style = {thick, dotted},
        ]
        coordinates {
          (67200,.75) (1075200,30.2)
        };
        \addlegendentry{$\mathcal O(n^{4/3})$};

      \end{loglogaxis}
    \end{tikzpicture}
    \caption{Solve time.}
  \end{subfigure}\\
  \vspace{\baselineskip}
  \begin{subfigure}[t]{.45\linewidth}
    \centering
    \begin{tikzpicture}[scale=0.8]
      \begin{loglogaxis} [
        legend pos = south east,
        xlabel = {$n$},
        ylabel = {$t_{q}$ (sec)}
        ]
        \addplot[
        color = black,
        mark = square,
        style = thick,
        ]
        coordinates {
          (67200,130.2) (151200, 260.8) (268800,470.1) (604800,1077.3) 
          (1075200,1936.4)
        };
        \addlegendentry{$n \propto k^2$};

        \addplot[
        color = gray,
        style = {thick, dashed},
        ]
        coordinates {
          (67200,130) (1075200,2080)
        };
        \addlegendentry{$\mathcal O(n)$};

        \addplot[
        color = gray,
        style = {thick, dotted},
        ]
        coordinates {
          (67200,130) (1075200,5239)
        };
        \addlegendentry{$\mathcal O(n^{4/3})$};

      \end{loglogaxis}
    \end{tikzpicture}
    \caption{Initialization time (quadrature generation).}
  \end{subfigure} \quad
  \begin{subfigure}[t]{.45\linewidth}
    \centering
    \begin{tikzpicture}[scale=0.8]
      \begin{loglogaxis} [
        legend pos = south east,
        xlabel = {$n$},
        ylabel = {$m_{f}$ (GB)}
        ]
        \addplot[
        color = black,
        mark = square,
        style = thick,
        ]
        coordinates {
          (67200,5.3) (151200, 10.8) (268800,18.2) (604800,38.1) 
          (1075200,65.5)
        };
        \addlegendentry{$n \propto k^2$};

        \addplot[
        color = gray,
        style = {thick, dashed},
        ]
        coordinates {
          (67200,5) (1075200,80)
        };
        \addlegendentry{$\mathcal O(n)$};

        \addplot[
        color = gray,
        style = {thick, dotted},
        ]
        coordinates {
          (67200,5) (1075200,201.5)
        };
        \addlegendentry{$\mathcal O(n^{4/3})$};

      \end{loglogaxis}
    \end{tikzpicture}
    \caption{Memory requirements.}
  \end{subfigure}
  \caption{Timings and memory usage for the FMM-LU solver for increasing
  frequency, based on data in Table~\ref{tab:freq}. The wavenumber~$k$ scales
  quadratically with~$n$, thereby keeping the approximate points-per-wavelength
  along the surface constant.}
  \label{fig:graphsPPW}
\end{figure}
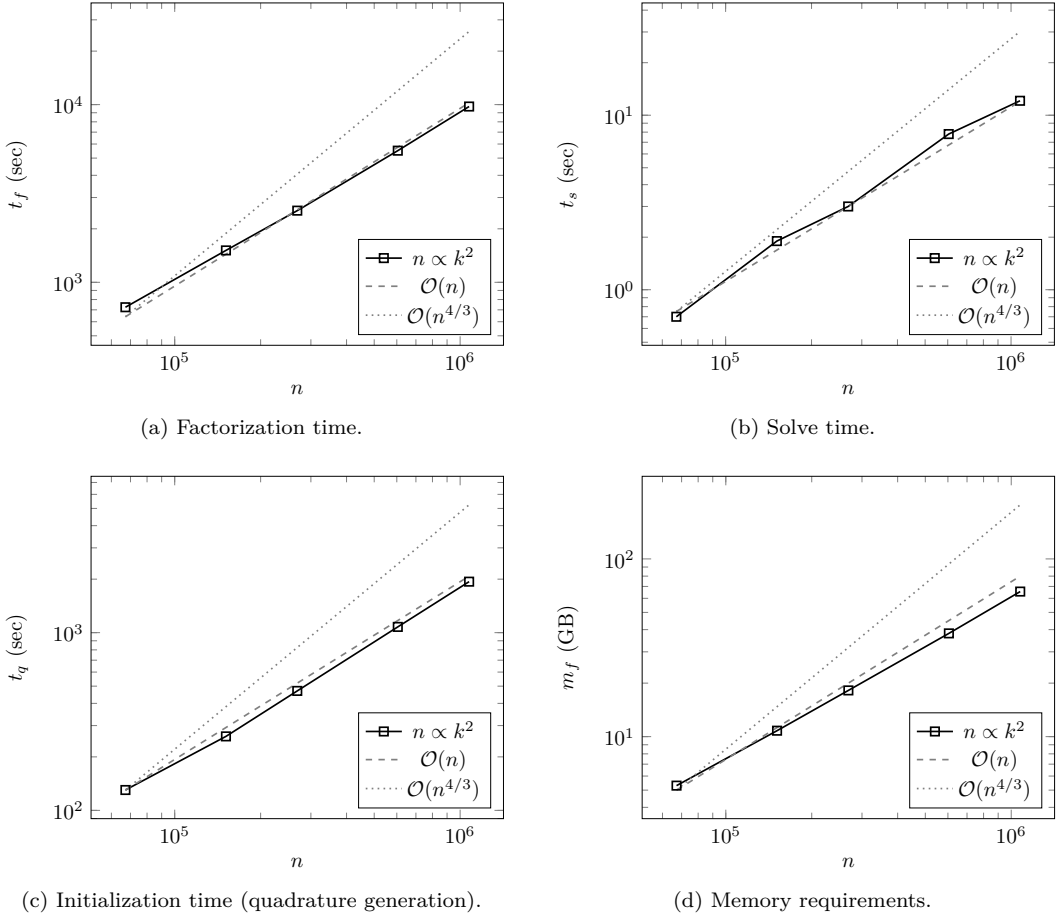
\clearpage
}

\subsection{Oscillatory problems sampled at fixed points per
  wavelength}
\label{subsec:num-freq}

In this section, we perform the accuracy test as in~\cref{subsec:num-conv} when
the surface $\Gamma$ is sampled at a fixed number of points per wavelength. For
these examples, $p=6$, $\varepsilon=5\times 10^{-7}$ and and
$k/\sqrt{\Npatches}$ is held fixed, so that $\ppwmin \approx 10$.
In~\cref{tab:freq}, we tabulate $\varepsilon_{a}$, $n_{0}$, $t_{f}$, $t_{s}$,
$t_{q}$, and $m_{f}$. The data are plotted in Figure~\ref{fig:graphsPPW}.
The error $\varepsilon_{a}$ remains nearly constant as
$k,\Npatches$ are increased since the surface is being sampled at nearly the
same number of points per wavelength. The size,~$n_{0}$, of the matrix at the root level
grows like $n^{4/9}$. Since the factorization time will typically be
dominated by the direct inversion at the root level as $k\to \infty$, the
computational complexity of constructing the factorization will be
$\cO(n_{0}^{3}) = \cO(n^{4/3})$. In the current experiments, we still
empirically observe an $\cO(n)$ scaling in the factorization time as the
frequency $k$ is not sufficiently large. The solve time, time for generating the
locally corrected quadrature, and the memory used continue to scale as $\cO(n)$
even in this setting. This can be explained by the fact that at the root level,
the solve time and memory requirement grow as $\cO(n_{0}^2) = \cO(n^{8/9})$, so
the dominant cost for both the solve time and memory used is at the leaf level
of the tree which is $\cO(n)$. Since we are using locally-corrected quadratures,
the number of near interactions that need to be precomputed are more or less
constant as $k,\Npatches \to \infty$, and thus the cost of computing them scales
like $\cO(n)$.

\subsection{Computing an azimuthal sonar cross section}
\label{sec:plane}

In this section, we compute the azimuthal monostatic sonar cross section for
sound soft acoustic scatterers. As before, let $\Gamma$ denote the boundary of
the scatterer whose interior is defined by $\Omega$. Let $\phi$ denote the
azimuthal angle in the $xy$ plane of the object, and let the azimuthal sonar
cross section in the $\phi$ direction be denoted by~$R(\phi)$. Furthermore,
let~$\mathcal{F}[u](\theta,\phi)$ denote the far field
pattern~\cite{colton_kress_inverse} of an outgoing solution~$u$ to the Helmholtz
equation:
\begin{equation}
  u(\bx) = \frac{e^{{ikr}}}{r} \cF [u] (\theta,\phi) + \cO\lp \frac{1}{r^{2}} \rp.
\end{equation}
Next, let~$u_{\phi_{0}}$ denote the outgoing scattered field due an incident plane wave
propagating in the direction~\mbox{$\bm{d} = (\cos(\phi_{0}),\sin{(\phi_{0})},0)$}.
Then, the monostatic cross section is given by
\begin{equation}
  R(\phi_0) = \cF[u_{\phi_{0}}](0,-\phi_0).
\end{equation}
That is to say: the monostatic cross section is the reflected far field
signature of the solution in the direction from which the incident plane wave
originated.

As before, we represent the scattered field~$u_{\phi_{0}}$ using the combined
field layer potential with unknown density~$\sigma_{\phi_{0}}$. Then, along the
boundary~$\Gamma$, the density~$\sigma_{\phi_{0}}(\bx)$ satisfies
\begin{equation}
  \frac{1}{2} \sigma_{\phi_{0}}(\bx) + \int_{\Gamma} K(\bx,\by) \, \sigma_{\phi_{0}}(\by) \, da(\by) = -e^{ik \bx \cdot \bm{d}},
\end{equation}
where the right hand side above is a plane wave propagating in the
direction~$\bm{d}$. It can then be shown that
the azimuthal sonar cross section~$R(\phi_{0})$ can be
expressed in terms of the solution~$\sigma_{\phi_{0}}$ as
\begin{equation}
    R(\phi_{0}) = \frac{-ik}{4\pi} \int_{\Gamma} e^{-i k \bx \cdot
      \bm{d}} \left( 1 - \bn(\bx) \cdot \bm{d} \right)
    \sigma_{\phi_{0}}(\bx) da(\bx),
\end{equation}
from which one can determine~$\mathcal{F}[u_{\phi_{0}}](0,-\phi_0)$.

In this example, we compute the azimuthal sonar cross section~$R$ of a model
airplane. The model airplane is~$49.3$ wavelengths long with a wingspan
of~$49.2$ wavelengths and a vertical height of~$13.7$ wavelengths. The plane
also has several multiscale features: 2 antennae on the top of the fuselage, and
1 \emph{control unit} on the bottom of the fuselage, see Figure~\ref{fig:a380}.
The ratio of the diameter of the largest triangle in the mesh to the smallest
triangle is~$\cO(10^3)$. The plane is discretized with~$\Npatches=125,344$ and~$p=4$,
resulting in~$n=1,253,440$ discretization points.

In Figure~\ref{fig:a380}, we plot the solution $\sigma_{\phi_{0}}$ for
$\phi_{0} = \pi/2$, and in Figure~\ref{fig:rad_600}, we plot the azimuthal sonar
cross section computed at $1000$ equispaced azimuthal angles on $(0,2\pi]$. We
also tested the accuracy of the factorization using the approach discussed
in~\cref{subsec:num-conv}. To generate the validation boundary data, we use $120$ point
sources located in the interior of the plane, and computed the accuracy of the
numerical solution on a~$101 \times 101$ lattice of targets on a slice which cuts
through the wing edge, whose normal is given by $(0,0,1)$. In this example, the
requested precision was set to $\varepsilon = 5\times 10^{-5}$, the
factorization was computed on a $40$ core linux workstation with $768$ GB RAM,
$t_{f} = 17810s$, $t_{s} = 32.2s$, $t_{q} = 60.5s$, $m_{f} = 459.5$GB,
$n_{0} = 17466$, and $\varepsilon_{a} = 7.4 \times 10^{-4}$.

\afterpage{

\begin{figure}[t]
  \centering
    \includegraphics[width=.95\linewidth]{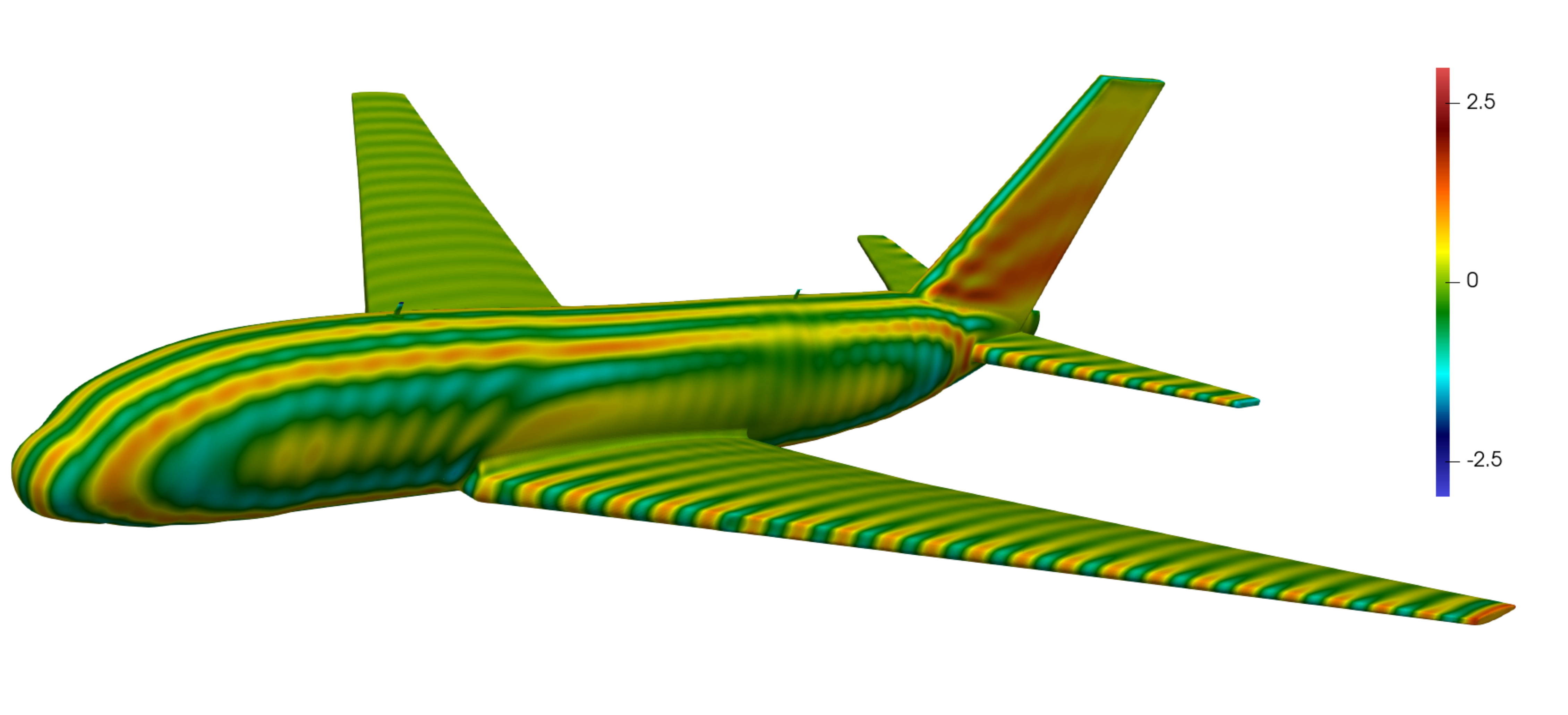}
  \caption{A multiscale model A380 aircraft. Plotted on surface is $\textrm{Re}(\sigma_{\pi/2})$}
  \label{fig:a380}
\end{figure}

\begin{figure}[b]
    \centering
    \includegraphics[width=0.65\linewidth]{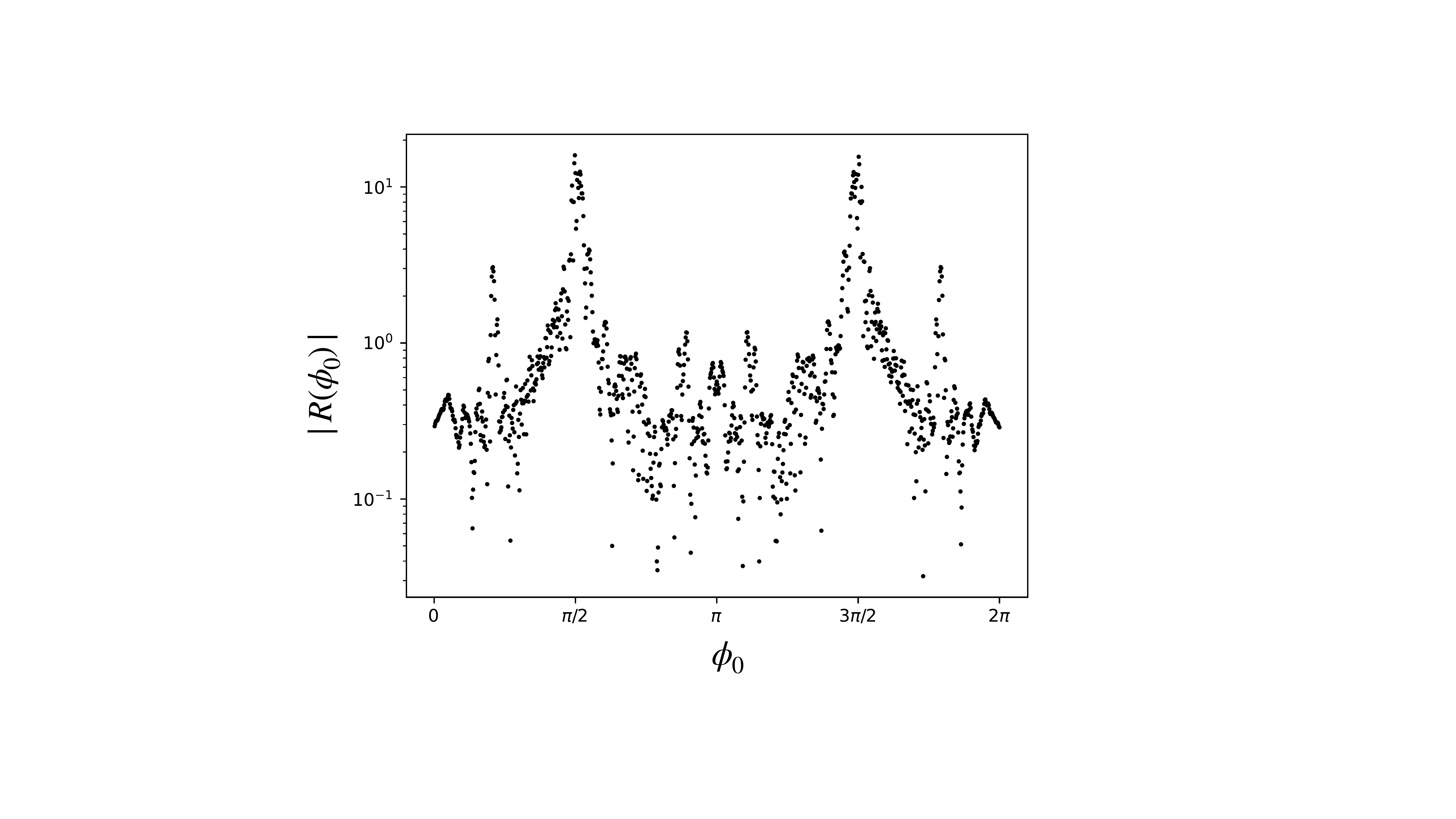}
    \caption{Azimuthal sonar cross section of the mock-up A380.}
    \label{fig:rad_600}
\end{figure}

\clearpage
}

\section{Conclusions} \label{sec:conclusions}


Fast direct solvers for boundary integral equations on complex surfaces in three
dimensions are just now beginning to be competitive with iterative techniques
coupled with fast-multipole methods (or related acceleration techniques). A
prototype for these schemes is the recursive skeletonization method introduced
by Martinsson and Rokhlin \cite{martinsson-2005}, which motivated the subsequent
hierarchical methods of
\cite{gillman-dir_hss-2012,ho-2012,kong-hodlrdir-2011,Mar-hss-2011} and is
related to the ${\cal H}$-matrix approach of
\cite{hackbusch-h-1999,hackbusch-h-2000,Borm-h-2003}. These methods are more or
less optimal in two dimensions but not in three dimensions because of the
high-rank interaction between adjacent surface patches. The methods of
\cite{ambikasaran-2014,Borm-h2-2010,hackbusch-h2-2000,
  Jiao-h2-ce-2017,minden-str_scel-2017,coulier-ifmm_prec-2017} which rely on
strong skeletonization (strong admissibility, mosaic-skeletonization,
${\cal H}^2$-matrix compression, or the inverse FMM) avoid direct compression of
these interactions, introducing more complex linear algebraic structures, with
one approach described in the body of the present paper.

While we have focused on the FMM-LU solver as an extension of the RS-S
method~\cite{minden-str_scel-2017} while using high-order accurate quadratures
at high precision, it should be noted that one can instead use a direct solver
at \emph{low precision} as a preconditioner for an iterative method. As noted
in~\cite{minden-str_scel-2017}, this may be a more practical approach because of
the reduction in memory requirements. At present, it seems as though the strong
skeletonization framework leads to methods with both the best asymptotic
scaling and the lowest memory requirements.

We have used the Helmholtz equation as our model problem in this work, and have
limited our attention to the low or moderate frequency regime. For highly
oscillatory problems, even well-separated blocks are of high rank and the
interpolative decomposition does not lead to significant compression.
Alternative compression techniques such as low rank directional
compression~\cite{borm_2016,borm_2015,engquist_2009} or butterfly
compression~\cite{guo_2016,liu_2016,li-2015} techniques are promising in this
regime.

There are a number of open problems in this rapidly evolving field: a better
understanding of the rank structure of the Schur complements that play a role in
the recursive solver, optimization of proxy surface selection, extension to the
high-frequency regime, and implementation on high-performance computing
architectures. Finally, a number of integral equation formulations of scattering
problems involve the {\em composition} of integral operators (see, for example,
the regularized combined field equation for sound-hard scatterers
\cite{bruno_2012,greengard2021fast}). In that setting, the entries of the system
matrix are not directly available in~$\cO(1)$ time. We are currently working on
a number of these issues.

\section*{Acknowledgments}

The authors would like to thank Felipe Vico for many useful conversations, and
James Bremer and Zydrunas Gimbutas for providing generalized Gaussian
quadratures that were used in the surface discretizations.

\bibliographystyle{abbrv}
\bibliography{lib}

\begin{thebibliography}{10}

\bibitem{ambikasaran2013n}
S.~Ambikasaran and E.~Darve.
\newblock {An $\mathcal{O}(N\log N)$ Fast Direct Solver for Partial
  Hierarchically Semi-Separable Matrices}.
\newblock {\em J. Sci. Comput.}, 57(3):477--501, 2013.

\bibitem{ambikasaran-2014}
S.~Ambikasaran and E.~Darve.
\newblock {The Inverse Fast Multipole Method}.
\newblock {\em arXiv preprint arXiv:1407.1572}, 2014.

\bibitem{ambikasaran_2016}
S.~Ambikasaran, D.~Foreman-Mackey, L.~Greengard, D.~W. Hogg, and M.~O'Neil.
\newblock {Fast Direct Methods for Gaussian Processes}.
\newblock {\em IEEE Trans. Pattern Anal. Mach. Intell.}, 38:252--265, 2016.

\bibitem{ambikasaran2016symmetric}
S.~Ambikasaran, M.~O'Neil, and K.~R. Singh.
\newblock Fast symmetric factorization of hierarchical matrices with
  applications.
\newblock {\em arXiv preprint arXiv:1405.0223}, 2016.

\bibitem{atkinson_1997}
K.~E. Atkinson.
\newblock {\em {The Numerical Solution of Integral Equations of the Second
  Kind}}.
\newblock Cambridge University Press, New York, NY, 1997.

\bibitem{Atkinson95}
K.~E. Atkinson and D.~Chien.
\newblock Piecewise polynomial collocation for boundary integral equations.
\newblock {\em SIAM J. Sci. Comput.}, 16:651--681, 1995.

\bibitem{shortcourse}
R.~Beatson and L.~Greengard.
\newblock A short course on fast multipole methods.
\newblock In M.~A. et~al., editor, {\em Wavelets, Multilevel Methods, and
  Elliptic PDEs}, pages 1--37. Oxford University Press, 1997.

\bibitem{Beb-hlu-2005}
M.~Bebendorf.
\newblock Hierarchical {LU} decomposition-based preconditioners for {BEM}.
\newblock {\em Computing}, 74(3):225--247, 2005.

\bibitem{Hakb-h2-lib}
S.~B{\"o}rm.
\newblock $\mathcal{H}$2lib package.
\newblock http://www.h2lib.org/.

\bibitem{Borm-h2-2010}
S.~B{\"o}rm.
\newblock {\em Efficient numerical methods for non-local operators:
  $\mathcal{H}^2$-matrix compression, algorithms and analysis}, volume~14.
\newblock European Mathematical Society, 2010.

\bibitem{borm_2015}
S.~B\"orm.
\newblock Directional ${H}^2$-matrix compression for high-frequency problems.
\newblock {\em Num. Lin. Alg. Appl.}, 24(6):e2112, 2017.

\bibitem{Borm-h-2003}
S.~B{\"o}rm, L.~Grasedyck, and W.~Hackbusch.
\newblock {Introduction to hierarchical matrices with applications}.
\newblock {\em Eng. Anal. Bound Elem.}, 27(5):405--422, 2003.

\bibitem{borm_2016}
S.~B{\"o}rm and J.~M. Melenk.
\newblock Approximation of the high-frequency {H}elmholtz kernel by nested
  directional interpolation: error analysis.
\newblock {\em Num. Math.}, 137(1):1--34, 2017.

\bibitem{borm-h2lu-2013}
S.~B{\"o}rm and K.~Reimer.
\newblock Efficient arithmetic operations for rank-structured matrices based on
  hierarchical low-rank updates.
\newblock {\em Comput. Vis. Sci.}, 16(6):247--258, 2013.

\bibitem{bremer_2012}
J.~Bremer.
\newblock On the {N}ystr\"om discretization of integral equations on planar
  curves with corners.
\newblock {\em Appl. Comput. Harm. Anal.}, 32:45--64, 2012.

\bibitem{bremer-2015}
J.~Bremer, A.~Gillman, and P.-G. Martinsson.
\newblock A high-order accelerated direct solver for integral equations on
  curved surfaces.
\newblock {\em BIT Num. Math.}, 55:367--397, 2015.

\bibitem{bremer_2012c}
J.~Bremer and Z.~Gimbutas.
\newblock A {N}ystr\"om method for weakly singular integral operators on
  surfaces.
\newblock {\em J. Comput. Phys.}, 231:4885--4903, 2012.

\bibitem{bremer_2013}
J.~Bremer and Z.~Gimbutas.
\newblock On the numerical evaluation of singular integrals of scattering
  theory.
\newblock {\em J. Comput. Phys.}, 251:327--343, 2013.

\bibitem{bremer}
J.~Bremer, Z.~Gimbutas, and V.~Rokhlin.
\newblock A nonlinear optimization procedure for generalized {G}aussian
  quadratures.
\newblock {\em SIAM J. Sci. Comput.}, 32(4):1761--1788, 2010.

\bibitem{bruno_2012}
O.~Bruno, T.~Elling, and C.~Turc.
\newblock Regularized integral equations and fast high-order solvers for
  sound-hard acoustic scattering problems.
\newblock {\em Int. J. Num. Meth. Engin.}, 91:1045--1072, 2012.

\bibitem{bruno_garza_2020}
O.~Bruno and E.~Garza.
\newblock A {C}hebyshev-based rectangular-polar integral solver for scattering
  by geometries described by non-overlapping patches.
\newblock {\em J. Comput. Phys.}, page 109740, 2020.

\bibitem{bruno2001fast}
O.~P. Bruno and L.~A. Kunyansky.
\newblock A fast, high-order algorithm for the solution of surface scattering
  problems: {B}asic implementation, tests, and applications.
\newblock {\em J. Comput. Phys.}, 169(1):80--110, 2001.

\bibitem{ChMing-dir_hss-2006}
S.~Chandrasekaran, M.~Gu, and T.~Pals.
\newblock A fast {ULV} decomposition solver for hierarchically semiseparable
  representations.
\newblock {\em SIAM J. Matrix Anal. A.}, 28(3):603--622, 2006.

\bibitem{chen-direct}
Y.~Chen.
\newblock {A fast, direct algorithm for the Lippmann-Schwinger integral
  equation in two dimensions}.
\newblock {\em Adv. Comput. Math.}, 16:175--190, 2002.

\bibitem{wideband3d}
H.~Cheng, W.~Y. Crutchfield, Z.~Gimbutas, L.~Greengard, J.~F. Ethridge,
  J.~Huang, V.~Rokhlin, N.~Yarvin, and J.~Zhao.
\newblock A wideband fast multipole method for the {H}elmholtz equation in
  three dimensions.
\newblock {\em J. Comput. Phys.}, 216:300--325, 2006.

\bibitem{cheng_2005}
H.~Cheng, Z.~Gimbutas, P.-G. Martinsson, and V.~Rokhlin.
\newblock On the compression of low rank matrices.
\newblock {\em SIAM J. Sci. Comput.}, 26(4):1389--1404, 2005.

\bibitem{gr-fmm_3d-1999}
H.~Cheng, L.~Greengard, and V.~Rokhlin.
\newblock A fast adaptive multipole algorithm in three dimensions.
\newblock {\em J. Comput. Phys.}, 155(2):468--498, 1999.

\bibitem{colton_kress}
D.~Colton and R.~Kress.
\newblock {\em Integral {E}quation {M}ethods in {S}cattering {T}heory}.
\newblock John Wiley \& Sons, Inc., 1983.

\bibitem{colton_kress_inverse}
D.~Colton and R.~Kress.
\newblock {\em Inverse {A}coustic and {E}lectromagnetic {S}cattering {T}heory}.
\newblock Springer, New York, NY, 2012.

\bibitem{CorMar-dir_hss-2014}
E.~Corona, P.-G. Martinsson, and D.~Zorin.
\newblock An $\mathcal{O}(n)$ direct solver for integral equations on the
  plane.
\newblock {\em Appl. Comput. Harmon. A.}, 2014.

\bibitem{coulier-ifmm_prec-2017}
P.~Coulier, H.~Pouransari, and E.~Darve.
\newblock {The inverse fast multipole method: Using a fast approximate direct
  solver as a preconditioner for dense linear systems}.
\newblock {\em SIAM J. Sci. Comput.}, 39:A761--A796, 2017.

\bibitem{treebook}
M.~de~Berg, O.~Cheong, M.~Kreveld, and M.~Overmars.
\newblock {\em Computational Geometry: Algorithms and Applications}.
\newblock Springer, 2008.

\bibitem{engquist_2009}
B.~Engquist and L.~Ying.
\newblock A fast directional algorithm for high frequency acoustic scattering
  in two dimensions.
\newblock {\em Commun. Math. Sci.}, 7:327--345, 2009.

\bibitem{erichsen1998quadrature}
S.~Erichsen and S.~A. Sauter.
\newblock {Efficient automatic quadrature in 3-D {G}alerkin {BEM}}.
\newblock {\em Comput. Methods Appl. Mech. Engrg.}, 157:215--224, 1998.

\bibitem{hps-direct}
A.~Gillman, A.~H. Barnett, and P.-G. Martinsson.
\newblock A spectrally accurate direct solution technique for frequency-domain
  scattering problems with variable media.
\newblock {\em BIT Num. Math.}, 55:141--170, 2015.

\bibitem{gillman-dir_hss-2012}
A.~Gillman, P.~M. Young, and P.-G. Martinsson.
\newblock {A direct solver with $O(N)$ complexity for integral equations on
  one-dimensional domains}.
\newblock {\em Frontiers of Mathematics in China}, 7(2):217--247, 2012.

\bibitem{gopal2020ls}
A.~Gopal and P.-G. Martinsson.
\newblock {An accelerated, high-order accurate direct solver for the
  Lippmann-Schwinger equation for acoustic scattering in the plane}.
\newblock {\em Adv. Comput. Math.}, 48, 2022.

\bibitem{greengard-2009}
L.~Greengard, D.~Gueyffier, P.-G. Martinsson, and V.~Rokhlin.
\newblock {Fast direct solvers for integral equations in complex
  three-dimensional domains}.
\newblock {\em Acta Numerica}, 18:243--275, 2009.

\bibitem{greengard2021fast}
L.~Greengard, M.~O'Neil, M.~Rachh, and F.~Vico.
\newblock Fast multipole methods for the evaluation of layer potentials with
  locally-corrected quadratures.
\newblock {\em J. Comput. Phys.: X}, 10:100092, 2021.

\bibitem{GrRo-fmm-1987}
L.~Greengard and V.~Rokhlin.
\newblock A fast algorithm for particle simulations.
\newblock {\em J. Comput. Phys.}, 73(2):325--348, Dec. 1987.

\bibitem{GrRo-fmm-1988}
L.~Greengard and V.~Rokhlin.
\newblock The rapid evaluation of potential fields in three dimensions.
\newblock In {\em Vortex methods}, pages 121--141. Springer, 1988.

\bibitem{gr-bvp}
L.~Greengard and V.~Rokhlin.
\newblock On the numerical solution of two-point boundary value problems.
\newblock {\em Comm. Pure Appl. Math.}, 44:419--452, 1991.

\bibitem{greengard-1997}
L.~Greengard and V.~Rokhlin.
\newblock {A new version of the Fast Multipole Method for the Laplace equation
  in three dimensions}.
\newblock {\em Acta Numerica}, 6:229--269, 1997.

\bibitem{gu1996efficient}
M.~Gu and S.~C. Eisenstat.
\newblock Efficient algorithms for computing a strong rank-revealing {QR}
  factorization.
\newblock {\em SIAM J. Sci. Comput.}, 17(4):848--869, 1996.

\bibitem{guo_2016}
H.~Guo, Y.~Liu, J.~Hu, and E.~Michielssen.
\newblock {A Butterfly-Based Direct Integral Equation Solver Using Hierarchical
  LU Factorization for Analyzing Scattering from Electrically Large Conducting
  Objects}.
\newblock {\em IEEE Trans. Antennas Propag.}, 65(9):4742--4750, 2017.

\bibitem{hackbusch-h-1999}
W.~Hackbusch.
\newblock A sparse matrix arithmetic based on $\mathcal{H}$-matrices. {Part I}:
  Introduction to $\mathcal{H}$-matrices.
\newblock {\em Computing}, 62(2):89--108, 1999.

\bibitem{hackbusch-h2-2000}
W.~Hackbusch, B.~Khoromskij, and S.~Sauter.
\newblock On $\mathcal{H}^2$-matrices.
\newblock In {\em H.-J. Bungartz, et al. (eds.), Lectures on Applied
  Mathematics}, pages 9--30. Springer-Verlag, Berlin Heidelberg, 2000.

\bibitem{hackbusch-h-2000}
W.~Hackbusch and B.~N. Khoromskij.
\newblock {A Sparse $\mathcal{H}$-Matrix Arithmetic}.
\newblock {\em Computing}, 64(1):21--47, 2000.

\bibitem{ho-2012}
K.~L. Ho and L.~Greengard.
\newblock {A fast direct solver for structured linear systems by recursive
  skeletonization}.
\newblock {\em SIAM J. Sci. Comput.}, 34:A2507--A2532, 2012.

\bibitem{ho-2014}
K.~L. Ho and L.~Ying.
\newblock Hierarchical interpolative factorization for elliptic operators:
  Integral equations.
\newblock {\em Comm. Pure Appl. Math.}, 69(7):1314--1353, 2016.

\bibitem{jiang2022skel}
M.~Jiang, Z.~Rong, X.~Yang, L.~Lei, P.~Li, Y.~Chen, and J.~Hu.
\newblock {Analysis of Electromagnetic Scattering From Homogeneous Penetrable
  Objects by a Strong Skeletonization-Based Fast Direct Solver}.
\newblock {\em IEEE Trans. Antennas Propag.}, 70(8):6883--6892, 2022.

\bibitem{kong-hodlrdir-2011}
W.~Y. Kong, J.~Bremer, and V.~Rokhlin.
\newblock An adaptive fast direct solver for boundary integral equations in two
  dimensions.
\newblock {\em Appl. Comput. Harm. Anal.}, 31(3):346--369, 2011.

\bibitem{kress_2014}
R.~Kress.
\newblock {\em {Linear Integral Equations}}.
\newblock Springer, New York, NY, 2014.

\bibitem{li-2015}
Y.~Li, H.~Yang, E.~Martin, K.~L. Ho, and L.~Ying.
\newblock Butterfly factorization.
\newblock {\em Multiscale Model. Simul.}, 13(2):714--732, 2015.

\bibitem{liberty2007randomized}
E.~Liberty, F.~Woolfe, P.-G. Martinsson, V.~Rokhlin, and M.~Tygert.
\newblock Randomized algorithms for the low-rank approximation of matrices.
\newblock {\em Proceedings of the National Academy of Sciences},
  104(51):20167--20172, 2007.

\bibitem{liubook}
Y.~Liu.
\newblock {\em Fast Multipole Boundary Element Method: Theory and Applications
  in Engineering}.
\newblock Cambridge University Press, 2008.

\bibitem{liu_2016}
Y.~Liu, H.~Guo, and E.~Michielssen.
\newblock {A HSS Matrix-Inspired Butterfly-Based Direct Solver for Analyzing
  Scattering from Two-dimensional Objects}.
\newblock {\em IEEE Antenn. Wirel. Pr.}, 16:1179--1183, 2016.

\bibitem{Jiao-h2-ce-2017}
M.~Ma and D.~Jiao.
\newblock {Accuracy Directly Controlled Fast Direct Solution of General
  ${\mathcal{ H}}^{2}$ -Matrices and Its Application to Solving Electrodynamic
  Volume Integral Equations}.
\newblock {\em IEEE Trans. Micro. Theory Tech.}, 66(1):35--48, 2018.

\bibitem{malhotra19}
D.~Malhotra, A.~Cerfon, L.-M. Imbert-G\'{e}rard, and M.~O'Neil.
\newblock Taylor states in stellarators: {A} fast high-order boundary integral
  solver.
\newblock {\em J. Comput. Phys.}, 397:108791, 2019.

\bibitem{Mar-hss-2011}
P.~G. Martinsson.
\newblock A fast randomized algorithm for computing a hierarchically
  semiseparable representation of a matrix.
\newblock {\em SIAM J. Matrix Anal. A.}, 32(4):1251--1274, 2011.

\bibitem{martinsson-book}
P.-G. Martinsson.
\newblock {\em Fast Direct Solvers for Elliptic PDEs}.
\newblock SIAM, 2019.

\bibitem{martinsson-2005}
P.-G. Martinsson and V.~Rokhlin.
\newblock {A fast direct solver for boundary integral equations in two
  dimensions}.
\newblock {\em J. Comput. Phys}, 205:1--23, 2005.

\bibitem{martinsson-hss_integr-2005}
P.-G. Martinsson and V.~Rokhlin.
\newblock A fast direct solver for boundary integral equations in two
  dimensions.
\newblock {\em J. Comput. Phys.}, 205(1):1--23, 2005.

\bibitem{martinsson2014}
P.-G. Martinsson, V.~Rokhlin, Y.~Shkolnisky, and M.~Tygert.
\newblock {\em {ID: A software package for low-rank approximation of matrices
  via interpolative decompositions, Version 0.4}}.
\newblock 0.4 edition, 2014.

\bibitem{minden-2016}
V.~Minden, A.~Damle, K.~L. Ho, and L.~Ying.
\newblock A technique for updating hierarchical skeletonization-based
  factorizations of integral operators.
\newblock {\em Multiscale Model. Simul.}, 14:42--64, 2016.

\bibitem{minden2017gp}
V.~Minden, A.~Damle, K.~L. Ho, and L.~Ying.
\newblock {Fast Spatial Gaussian Process Maximum Likelihood Estimation via
  Skeletonization Factorizations}.
\newblock {\em Multiscale Model. Simul.}, 15(4):1584--1611, 2017.

\bibitem{minden-str_scel-2017}
V.~Minden, K.~L. Ho, A.~Damle, and L.~Ying.
\newblock A recursive skeletonization factorization based on strong
  admissibility.
\newblock {\em Multiscale Model. Simul.}, 15(2):768--796, 2017.

\bibitem{Ro-fmm-1985}
V.~Rokhlin.
\newblock Rapid solution of integral equations of classical potential theory.
\newblock {\em J. Comput. Phys.}, 60(2):187--207, 1985.

\bibitem{Siegel2018ALT}
M.~Siegel and A.-K. Tornberg.
\newblock A local target specific quadrature by expansion method for evaluation
  of layer potentials in 3d.
\newblock {\em J. Comput. Phys.}, 364:365--392, 2018.

\bibitem{solovyev-hss-2014}
S.~Solovyev.
\newblock {Multifrontal Hierarchically Solver for 3D Discretized Elliptic
  Equations}.
\newblock In {\em International Conference on Finite Difference Methods}, pages
  371--378. Springer, 2014.

\bibitem{sr-bvp}
P.~Starr and V.~Rokhlin.
\newblock {On the numerical solution of two-point boundary value problems II}.
\newblock {\em Comm. Pure Appl. Math.}, 47:1117--1159, 1994.

\bibitem{sushnikova-ce-2018}
D.~A. Sushnikova and I.~V. Oseledets.
\newblock "{C}ompress and eliminate" solver for symmetric positive definite
  sparse matrices.
\newblock {\em SIAM J. Sci. Comput.}, 40(2):A1742--A1762, 2018.

\bibitem{vioreanu_2014}
B.~Vioreanu and V.~Rokhlin.
\newblock {Spectra of Multiplication Operators as a Numerical Tool}.
\newblock {\em SIAM J. Sci. Comput.}, 36:A267--A288, 2014.

\bibitem{Wala2018}
M.~Wala and A.~Kl{\"o}ckner.
\newblock A fast algorithm for quadrature by expansion in three dimensions.
\newblock {\em J. Comput. Phys.}, 388:655--689, 2018.

\bibitem{Wala2020}
M.~Wala and A.~Kl{\"o}ckner.
\newblock Optimization of fast algorithms for global quadrature by expansion
  using target-specific expansions.
\newblock {\em J. Comput. Phys.}, 403, 2020.

\bibitem{wu2020corrected}
B.~Wu and P.-G. Martinsson.
\newblock {Corrected Trapezoidal Rules for Boundary Integral Equations in Three
  Dimensions}.
\newblock {\em Num. Math.}, 149:1025--1071, 2021.

\bibitem{XiaSh-hss-2009}
J.~Xia, S.~Chandrasekaran, M.~Gu, and X.~S. Li.
\newblock Superfast multifrontal method for large structured linear systems of
  equations.
\newblock {\em SIAM J. Matrix Anal. A.}, 31(3):1382--1411, 2009.

\bibitem{ChLi-hss-2007}
J.~Xia, S.~Chandrasekaran, M.~Gu, and X.~S. Li.
\newblock Fast algorithms for hierarchically semiseparable matrices.
\newblock {\em Numer. Linear Algebr.}, 17(6):953--976, 2010.

\bibitem{xing2020proxy}
X.~Xing and E.~Chow.
\newblock {Interpolative Decomposition via Proxy Points for Kernel Matrices}.
\newblock {\em SIAM J. Matrix Anal. Appl.}, 41(1):221--243, 2020.

\bibitem{ye2020proxy}
X.~Ye, J.~Xia, and L.~Ying.
\newblock {Analytical Low-Rank Compression via Proxy Point Selection}.
\newblock {\em SIAM J. Matrix Anal. Appl.}, 41(3):1059--1085, 2020.

\bibitem{ying}
L.~Ying, G.~Biros, and D.~Zorin.
\newblock A high-order 3{D} boundary integral equation solver for elliptic
  {PDE}s in smooth domains.
\newblock {\em J. Comput. Phys.}, 219(1):247--275, 2006.

\end{thebibliography}

\end{document}